\documentclass{amsart}
\title%
[Values of symmetric power $L$-functions and Hecke eigenvalues]%
{Special values of symmetric power $L$-functions and Hecke eigenvalues}
\date{\today}
\usepackage[british]{babel}
\usepackage[T1]{fontenc}
\usepackage{amssymb}
 \usepackage{url}
 \usepackage[dvips]{graphicx}
\theoremstyle{plain}
\newtheorem{propint}{Proposition}

\newtheorem{thmint}[propint]{Theorem}

\newtheorem{corint}[propint]{Corollary}

\newtheorem{lemint}[propint]{Lemma}

\newtheorem*{HypSym}{Hypothesis $\sym^m(N)$}
\newtheorem*{HypLS}{Hypothesis $LSZ^m(N)$}
\theoremstyle{remark}
\newtheorem*{remint}{Remark}
\theoremstyle{plain}
\newtheorem{cor}{Corollary}
\newtheorem{lem}[cor]{Lemma}
\newtheorem{prop}[cor]{Proposition}
\theoremstyle{remark}
\newtheorem{rem}[cor]{Remark}
\DeclareFontFamily{U}{bbold}{}
\DeclareFontShape{U}{bbold}{m}{n}{<6> <7> <8> <9> <10> <12> <14> <17> bbold10}{}
\DeclareSymbolFont{gras}{U}{bbold}{m}{n}
\newcommand*{\abs}[1]{\lvert #1\rvert}
\newcommand{\bigtimes}[2]{\mathop{\mathrm{X}}_{#1}^{#2}}
\newcommand*{\C}{\mathbb{C}}
\def\dd{\mathop{\mathrm{d}\null}%
       \mskip-\thinmuskip\mathord{\null}}
\renewcommand*{\epsilon}{\varepsilon}
\DeclareMathOperator{\Err}{Err}
\renewcommand*{\Gamma}{\varGamma}
\newcommand*{\gl}{\mathrm{GL}}
\let\plainIm = \Im
\def\Im{\mathop{\plainIm\mkern-2mu \mathit{m}}\nolimits}

\newcommand*{\LB}{X}
\newcommand*{\N}{\mathbb{Z}_{\geq 0}}
\newcommand*{\Ncri}[1]{\mathcal{N}\left(#1\right)}
\newcommand*{\omegam}{\omega^*}
\newcommand*{\pk}{\mathcal{H}}
\newcommand*{\pr}{\mathrm{H}^*_{2}(N)}
\newcommand*{\prgde}[1]{\mathrm{H}^{*+}_{2}(N;C,\sym^{#1})}
\newcommand*{\prpte}[1]{\mathrm{H}^{*-}_{2}(N;C,\sym^{#1})}
\newcommand*{\prpn}[2]{\mathrm{H}^*_{#1}(#2)}
\newcommand*{\prem}{\mathcal{P}}
\newcommand*{\prm}{\mathrm{H}^{-}_{m}(N;\eta)}
\newcommand*{\prpdd}[1]{\mathrm{H}^{-}_{#1}(N;\eta)}
\newcommand*{\prd}{\mathrm{H}^{-}_{2,m}(N;\eta)}
\newcommand*{\prp}{\mathrm{H}^{+}_{m}(N;\eta)}
\newcommand*{\R}{\mathbb{R}}
\let\plainRe= \Re
\renewcommand*{\Re}%
      {\mathop{\plainRe\mkern-2mu \mathit{e}}\nolimits}
\DeclareMathOperator{\st}{St}
\newcommand*{\su}{\mathrm{SU}(2)}
\DeclareMathOperator{\sym}{Sym}
\DeclareMathOperator{\tr}{tr}
\DeclareMathSymbol{\un}{\mathord}{gras}{'061}
\def\vec#1{\ensuremath{\mathchoice
           {\mbox{\boldmath$\displaystyle\mathbf{#1}$}}
           {\mbox{\boldmath$\textstyle\mathbf{#1}$}}
           {\mbox{\boldmath$\scriptstyle\mathbf{#1}$}}
           {\mbox{\boldmath$\scriptscriptstyle\mathbf{#1}$}}}}
\newcommand*{\Z}{\mathbb{Z}}
\renewcommand*{\Sigma}{\varSigma}
\author[Emmanuel {\sc Royer}]{Emmanuel {\sc Royer}}
\address{%
Emmanuel {\sc Royer}\\
Laboratoire de math\'ematiques,\\ 
UMR6620 UBP~CNRS, \\
Universit\'e Blaise Pascal,\\
 Campus universitaire des C\'ezeaux, \\
F--63177 Aubi\`ere Cedex, \\
France%
}
\email{{emmanuel.royer@math.univ-bpclermont.fr}}
\urladdr{http://carva.org/emmanuel.royer}
\author[Jie {\sc Wu}]{Jie {\sc Wu}}
\address{%
Jie {\sc Wu}\\
Institut \'Elie Cartan, \\
UMR7502 UHP~CNRS~INRIA,\\%
 Universit\'e Henri Poincar\'e, Nancy 1, \\
F--54506 Vand{\oe}uvre-l\'es-Nancy, \\
France}
\email{{wujie@iecn.u-nancy.fr}}

\begin{document}

\begin{abstract}
We compute the moments of $L$-functions of symmetric powers of modular forms at the edge of the critical strip, twisted by the 
central value of the $L$-functions of modular forms. We show that, in the case of even powers, it is equivalent to twist by the value at 
the edge of the critical strip of the symmetric square $L$-functions. We deduce information on the size of symmetric power $L$-functions 
at the edge of the critical strip in subfamilies. In a second part, we study the distribution of small and large Hecke eigenvalues.
We deduce information on the simultaneous extremality conditions on the values of $L$-functions of symmetric powers of modular 
forms at the edge of the critical strip. 
\end{abstract}

\maketitle

\tableofcontents

\section{Introduction}

The values of $L$-functions at the edge of the critical strip have been extensively studied. The work on their
distributions originates with Littlewood \cite{L28}. In the case of %
Dirichlet $L$-functions, his %
work has been extended by Elliott \cite{ell73} and more recently by Montgomery \& Vaughan \cite{MV99} and
Granville \& Soundararajan \cite{GS03}. In the case of symmetric square $L$-functions of modular forms, the first
results are due to Luo \cite{Luo99}, \cite{Luo00}. They have been developed by the first author \cite{roy01} and the
authors \cite{RW03} in the analytic aspect and by the first author \cite{roy03} and Habsieger \& the first author \cite{MR2139690}
in the combinatorial aspect. These developments have been recently widely extended by Cogdell \& Michel \cite{CoMi04} who
studied the distribution for all the symmetric power $L$-functions.

The values of $L$-functions of modular forms at the centre of the critical strip are much more difficult to catch. The difficulty of
the computation of their moments increases dramatically with the order of the moments (see, e.g., \cite{KMV00}) and these
moments are subject to important conjectures \cite{CFKRS03}, \cite{CFKRS05}. Good bounds for the size of these values have important
consequences. A beautiful one is the following, due to Iwaniec \& Sarnak \cite{IS00}. Denote by
$\pr$\label{page:un} the set of primitive forms of weight $2$ over $\Gamma_0(N)$ (this is the Hecke eigenbasis of the space %
of parabolic newforms of weight 
$2$ over $\Gamma_0(N)$, normalised so that the first Fourier coefficient is one) and let $\epsilon_f(N)$ be the sign of the %
functional equation 
satisfied by 
the $L$-function, $L(s,f)$, of $f\in\pr$. Our $L$-functions are normalised so that $0\leq\Re s\leq 1$ is the critical strip. 
Then it is shown that
\[%
\liminf_{N\to\infty}
\frac{\#\left\{
f\in\pr \colon \epsilon_f(N)=1,\, L\left(\frac{1}{2},f\right)\geq(\log N)^{-2}
\right\}}{\#\{f\in\pr \colon \epsilon_f(N)=1\}}
\geq
c=\frac{1}{2}.
\]
If we could replace $c=1/2$ by $c>1/2$, then there would exist no Landau-Siegel zero for Dirichlet $L$-functions. It is expected
that one may even take $c=1$. The meaning of this expectation is that, if $L(1/2,f)\neq 0$ (which is not the case when $\epsilon_f(N)\neq 1$),
then $L(1/2,f)$ is not too small.

In this paper, we compute (see Theorem~\ref{thmint:A} and Proposition~\ref{propint:B}) the moments of symmetric power $L$-functions
at $1$ twisted by the value at $1/2$ of modular forms $L$-functions, that is
\begin{equation}\label{eq:momun}
\sum_{f\in\pr}\omegam(f)L\left(\frac{1}{2},f\right)L(1,\sym^mf)^z \quad (z\in\C)
\end{equation}
where $\omegam$ is the usual harmonic weight (see \eqref{eq:hawe}). Comparing (see Theorem~\ref{thmint:B} and Proposition~\ref{propint:D}) 
with the moments of symmetric power $L$-functions at $1$ twisted by the value of the symmetric square $L$-function at $1$, that is
\begin{equation}\label{eq:momde}
\sum_{f\in\pr}\omegam(f)L(1,\sym^2f)L(1,\sym^mf)^z  \quad (z\in\C),
\end{equation}
we show in Corollary~\ref{corint:identity} that \eqref{eq:momun} and \eqref{eq:momde} have asymptotically (up to a
multiplicative factor $1/\zeta(2)$) the same value when 
the power $m$ is even. This equality 
is astonishing since half of the values $L(1/2,f)$ are expected to be $0$ whereas $L(1,\sym^2f)$ is always positive. Since it
is even expected that $L(1,\sym^2f)\gg[\log\log(3N)]^{-1}$, it could suggest that $L(1/2,f)$ is large when not vanishing.

Our computations also yield results on the size of $L(1,\sym^mf)$ when subject to condition on the nonvanishing of $L(1/2,f)$
(see Corollary~\ref{lem:dblcdn}) or to extremality conditions for another symmetric power $L$-function (see Propositions~\ref{propint:simul}
and~\ref{propint:simulde}).

Before giving precisely the results, we introduce a few basic facts needed for the exposition.
More details shall be given in Section~\ref{sec:tools}.
Let $f$ be an element of the set $\pr$ of primitive forms of weight $2$
and squarefree level $N$ (i.e., over $\Gamma_0(N)$ and without nebentypus). It admits a
Fourier expansion 
\begin{equation}\label{eq:deff}
f(z)=: \sum_{n=1}^{+\infty}\lambda_f(n)\sqrt{n}e^{2\pi inz}
\end{equation}
in the upper half-plane $\pk$. Denote by $\st$ the standard representation of $\su$, 
\[%
\begin{array}{ccccc}
\st & \colon & \su & \to     & \gl(\C^2) \\
    &        & M & \mapsto & 
\begin{array}{ccc}
\C^2 & \to & \C^2\\
x    & \mapsto & Mx
\end{array}
\end{array}
\]
(for the basics on representations, see, e.g., \cite{vil68}).
If $\rho$ is a representation of $\su$ and $I$ is the identity matrix, define, 
for each $g\in\su$
\begin{equation}\label{eq:fato}
D(X,\rho,g):= \det[I-X\rho(g)]^{-1}.
\end{equation}
Denote by $\chi_\rho$\label{def:char} the character of $\rho$.
By Eichler \cite{eic54} and Igusa \cite{igu59}, we know that for every prime number $p$ not dividing the level, 
$\abs{\lambda_f(p)}\leq 2$ so that there exists $\theta_{f,p}\in[0,\pi]$ such that
\begin{equation*}
\lambda_f(p)=\chi_{\st}[g(\theta_{f,p})]
\end{equation*}
where 
\begin{equation}\label{eq:defgt}
g(\theta):= \begin{pmatrix}e^{i\theta} & 0\\0 & e^{-i\theta}\end{pmatrix}
\end{equation}
(in other words, $\lambda_f(p)=2\cos\theta_{f,p}$: this is the special case for weight $2$ forms of
the Ramanujan conjecture proved by Deligne for every weights).
Denote by $\prem$ the set of prime numbers.
Consider the symmetric power $L$-functions 
of $f$ defined for every integer $m\geq 0$ by
\begin{equation}\label{eq:defsym}
L(s,\sym^mf):= \prod_{p\in\prem}L_p(s,\sym^mf)
\end{equation}
where
\[%
L_p(s,\sym^mf):= D[p^{-s},\sym^m,g(\theta_{f,p})]
\]
if $p$ is coprime to the level $N$ and
\[%
L_p(s,\sym^mf):= [1-\lambda_f(p^m)p^{-s}]^{-1}
\]
otherwise. Here $\sym^m$ denotes the composition of the
$m$th symmetric power representation of $\gl(2)$ and the standard representation of $\su$. In
particular $\sym^0(g)=1$ for all $g\in\gl(2)$ so that $\sym^0$ is the trivial irreducible representation
and $L(s,\sym^0f)$ is the Riemann $\zeta$ function. 

We shall give all our results in a restrictive range for $m$. If we assume two standard 
hypothesis -- see Section~\ref{sec:tsh} --  the restriction is no longer necessary, i.e., all results
are valid for every integer $m\geq 1$.

\subsection{Twisted moments}

For each squarefree positive integer $N$, each positive integer $m$ and each complex number $z$, define
\begin{equation}\label{eq:dtm}
\LB_m^z(N):= \sum_{n=1}^{+\infty}\frac{\tau_z(n)}{n^{m/2+1}}\sum_{q=1}^{+\infty}\frac{\square_N(n^mq)}{q}
\end{equation}
where $\tau_z$ and $\square_N$ are defined by 
\begin{align}
\label{eq:zz}
\sum_{n=1}^{+\infty}\frac{\tau_z(n)}{n^s} &:=  \zeta(s)^z ,\\
\label{eq:caN}
\sum_{n=1}^{+\infty}\frac{\square_N(n)}{n^s}  &:=  \zeta_N(2s) 
:= \frac{\zeta(2s)}{\zeta^{(N)}(2s)}
:=  \prod_{\substack{p\in\prem\\ p\mid N}}
\left(1-\frac{1}{p^{2s}}\right)^{-1},
\end{align}
and
\begin{multline}\label{eq:eli1}
L^{1,z}\left(\frac{1}{2},1;\st,\sym^m;N\right)
\\
:= %
\LB_m^z(N)\prod_{\substack{p\in\prem\\ (p,N)=1}}
\int_{\su}D(p^{-1/2},\st,g)D(p^{-1},\sym^m,g)^z\dd g
\end{multline}
where $\dd g$ stands for the Haar measure on $\su$. In the special case $N=1$ write
\begin{multline}\label{eq:eli2}
L^{1,z}\left(\frac{1}{2},1;\st,\sym^m\right)
\\ :=%
\prod_{p\in\prem}%
\int_{\su}D(p^{-1/2},\st,g)D(p^{-1},\sym^m,g)^z\dd g.
\end{multline}
We also use the usual harmonic weight on the space of cuspidal forms
\begin{equation}\label{eq:hawe}
\omegam(f):= \frac{1}{4\pi(f,f)}\cdot\frac{N}{\varphi(N)}
\end{equation}
where $(f,f)$ is the Petersson norm of $f$ and $\varphi$ is Euler's %
totient function. We slightly change the usual definition to obtain
\[%
\lim_{N\to +\infty}\sum_{f\in\pr}\omegam(f)=1
\]
as $N$ runs over squarefree integers (see Lemma~\ref{lem:ils} with $m=n=1$) in order to
obtain an asymptotic average operator. We note $\log_n$ for the logarithm iterated $n$ times: 
$\log_1:= \log$ and $\log_{n+1}:= \log\circ\log_n$.
Our first result computes the twisted moments as in \eqref{eq:momun}.
\begin{thmint}\label{thmint:A}
Let $m\in\{1,2,4\}$. There exist two real numbers $c>0$ and $\delta>0$ such that, for
any squarefree integer $N\geq 1$, for any complex number $z$ verifying
\[%
\abs{z}\leq c\frac{\log(2N)}{\log_2(3N)\log_3(20N)}
\]
the following estimate holds:
\begin{multline*}
\sum_{f\in\pr}\omegam(f)L\left(\frac{1}{2},f\right)L(1,\sym^mf)^z
\\=
L^{1,z}\left(\frac{1}{2},1;\st,\sym^m;N\right)
+
O_m\left(\exp\left[-\delta\frac{\log(2N)}{\log_2(3N)}\right]\right)
\end{multline*}
with an implicit constant depending only on $m$.
\end{thmint}
Moreover, we obtain an asymptotic expression as $N$ tends to infinity in the next proposition. Define, for each
function $g\colon\Z_{>0}\to\R^+$, the set
\begin{equation}\label{eq:defng}
\Ncri{g}:=\{N\in\Z_{>0} \colon \mu(N)^2=1,\, P^-(N)\geq g(N)\}
\end{equation}
where $P^-(N)$\label{page:cin} is the smallest prime divisor of $N$ with the convention $P^-(1):=+\infty$, $\omega(N)$
is the number of distinct prime divisors of $N$ and $\mu$ is the M{\"o}bius function.
\begin{propint}\label{propint:B}
Let $\xi$ be a function such that $\xi(N)\to+\infty$ as $N\to+\infty$. Then
\[%
L^{1,z}\left(\frac{1}{2},1;\st,\sym^m;N\right)
=
L^{1,z}\left(\frac{1}{2},1;\st,\sym^m\right)
[1+o_m(1)]
\]
uniformly for
\begin{equation}\label{eq:j18}
\begin{cases}
N\in\Ncri{%
\xi(\cdot)
\max\bigl\{
\omega(\cdot),\,
[(\abs{z}+1)\omega(\cdot)]^{2/3},\, 
(\abs{z}+1)\omega(\cdot)^{1/2} 
\bigl\}
}
,\text{}\\
\abs{z}\leq  c \log(2N)/[\log_2(3N)\log_3(20N)].
\end{cases}
\end{equation}
\end{propint}
\begin{remint}
Condition \eqref{eq:j18} is certainly satisfied for
\[%
N\in\Ncri{\log^{3/2}}
\qquad\hbox{and}\qquad
\abs{z}\leq  c \log(2N)/[\log_2(3N)\log_3(20N)].
\]
\end{remint}
For a comparison of the behaviour of $L(1/2,f)$ and $L(1,\sym^2f)$ we next compute the moments of $L(1,\sym^mf)$
twisted by $L(1,\sym^2f)$. Define
\begin{equation}\label{eq:lenti}
\LB_{2,m}^{1,z}(N):= 
\zeta_N(2)
\sum_{n=1}^{+\infty}\frac{\tau_z(n)\square_N(n^m)}{n^{m/2+1}}
\end{equation}
and
\begin{multline}\label{eq:heure}
L^{1,z}\left(1,1;\sym^2,\sym^m;N\right)
\\
:= 
\LB_{2,m}^{1,z}(N)\prod_{\substack{p\in\prem\\ (p,N)=1}}
\int_{\su}D(p^{-1},\sym^2,g)D(p^{-1},\sym^m,g)^z\dd g.
\end{multline} 
For the special case $N=1$ we get
\begin{multline}\label{eq:charo}
L^{1,z}\left(1,1;\sym^2,\sym^m\right)
\\
:= %
\prod_{p\in\prem}%
\int_{\su}D(p^{-1},\sym^2,g)D(p^{-1},\sym^m,g)^z\dd g.
\end{multline}

\begin{thmint}\label{thmint:B}
Let $m\in\{1,2,4\}$. There exist two real numbers $c>0$ and $\delta>0$ such that, for
any squarefree integer $N\geq 1$, for
any complex number $z$ verifying
\[%
\abs{z}\leq c\frac{\log(2N)}{\log_2(3N)\log_3(20N)}
\]
the following estimate holds:
\begin{multline*}
\sum_{f\in\pr}\omegam(f)L\left(1,\sym^2f\right)L(1,\sym^mf)^z
\\=
L^{1,z}\left(1,1;\sym^2,\sym^m;N\right)
+
O\left(\exp\left[-\delta\frac{\log(2N)}{\log_2(3N)}\right]\right)
\end{multline*}
with an implicit constant depending only on $m$.
\end{thmint}
Again, we obtain an asymptotic expansion in the following proposition.
\begin{propint}\label{propint:D}
Let $\xi$ be a function such that $\xi(N)\to+\infty$ as $N\to+\infty$. Then
\[%
L^{1,z}\left(1,1;\sym^2,\sym^m;N\right)
=
L^{1,z}\left(1,1;\sym^2,\sym^m\right)
[1+o_m(1)]
\]
uniformly for
\begin{equation}\label{eq:j19}
\begin{cases}
N\in\Ncri{%
\xi(\cdot)
\max\bigl\{
\omega(\cdot)^{1/2},\, 
[(\abs{z}+1)\omega(\cdot)]^{2/(m+2)}
\bigl\}
}
,\text{}\\
\abs{z}\leq  c \log(2N)/[\log_2(3N)\log_3(20N)].
\end{cases}
\end{equation}
\end{propint}
\begin{remint}
Condition \eqref{eq:j19} is certainly satisfied for
\[%
N\in\Ncri{\log^{4/3}}
\qquad\hbox{and}\qquad
\abs{z}\leq  c \log(2N)/[\log_2(3N)\log_3(20N)].
\]
\end{remint}
From Theorems~\ref{thmint:A} and~\ref{thmint:B} and
\begin{multline*}
\prod_{p\in\prem}\int_{\su}D(p^{-1},\sym^{2m},g)^zD(p^{-1/2},\st,g)\dd g
\\=
\frac{1}{\zeta(2)}\prod_{p\in\prem}\int_{\su}D(p^{-1},\sym^{2m},g)^zD(p^{-1},\sym^2,g)\dd g
\end{multline*}
(see Lemma~\ref{lem:j46}), we deduce the following astonishing result.
\begin{corint}\label{corint:identity}
Let $m\in\{2,4\}$.
For any $N\in\Ncri{\log}$ and $f\in\pr$, for any $z\in\C$, the following estimate holds:
\begin{multline*}
\lim_{\substack{%
N\to\infty\\ N\in\Ncri{\log}
}}
\sum_{f\in\pr}\omegam(f)L\left(\frac{1}{2},f\right)L(1,\sym^{m}f)^z
\\=
\lim_{\substack{%
N\to\infty\\ N\in\Ncri{\log}
}}
\frac{1}{\zeta(2)}\sum_{f\in\pr}\omegam(f)L(1,\sym^2f)L(1,\sym^{m}f)^z.
\end{multline*}
\end{corint}
This identity is not valid when considering an odd symmetric power
of $f$. For example,
\begin{multline}\label{eq:su}
\lim_{\substack{%
N\to\infty\\ N\in\Ncri{\log}
}}
\sum_{f\in\pr}\omegam(f)L\left(\frac{1}{2},f\right)L(1,f)
\\
=%
\prod_{p\in\prem}\left(1+\frac{1}{p^{3/2}}+O\left(\frac{1}{p^2}\right)\right)
\end{multline}
and 
\begin{equation}\label{eq:sd}
\lim_{\substack{%
N\to\infty\\ N\in\Ncri{\log}
}}
\sum_{f\in\pr}\omegam(f)L\left(1,\sym^2f\right)L(1,f)
=
\prod_{p\in\prem}\left(1+O\left(\frac{1}{p^2}\right)\right)
\end{equation}
so that the quotient of \eqref{eq:su} by \eqref{eq:sd} is
\[%
\prod_{p\in\prem}\left(1+\frac{1}{p^{3/2}}+O\left(\frac{1}{p^2}\right)\right)
\]
whereas
\begin{equation}\label{eq:st}
\lim_{\substack{%
N\to\infty\\ N\in\Ncri{\log}
}}
\sum_{f\in\pr}\omegam(f)L\left(\frac{1}{2},f\right)L(1,\sym^{3}f)
=
\prod_{p\in\prem}\left(1+O\left(\frac{1}{p^2}\right)\right)
\end{equation}
and
\begin{multline}\label{eq:sq}
\lim_{\substack{%
N\to\infty\\ N\in\Ncri{\log}
}}
\frac{1}{\zeta(2)}\sum_{f\in\pr}\omegam(f)L(1,\sym^2f)L(1,\sym^{3}f)
\\
=%
\prod_{p\in\prem}\left(1+O\left(\frac{1}{p^2}\right)\right)
\end{multline}
so that the quotient of \eqref{eq:st} by \eqref{eq:sq} is
\[%
\prod_{p\in\prem}\left(1+O\left(\frac{1}{p^2}\right)\right).
\]
The key point of Corollary~\ref{corint:identity} is the fact that the coefficients appearing in the series expansion of
$D(X,\sym^{2m},g)$ have only even harmonics -- see equations \eqref{eq:j23} and \eqref{eq:j27}. See Remark~\ref{rem:plana} for further details. 

\subsection{Extremal values}

The size of the values $L(1,\sym^mf)$ in the family $\pr$ is now well studied after works of %
Goldfeld, Hoffstein \& Lieman \cite{ghl94}, Royer \& Wu \cite{RW03}, Cogdell \& Michel \cite{CoMi04} %
and Lau \& Wu \cite{LW07} (among others). The aim of Proposition~\ref{propint:asymp} and Corollary~\ref{lem:dblcdn} %
is to study the extremal values in some smaller family. More precisely, we study the extremal values in families %
determined by the nonvanishing of $L\left(\frac{1}{2},f\right)$ and show that the extremal values are the same than in %
the full family.

We begin in studying the asymptotic behaviour, as the order $z$ tends to $\pm\infty$ in %
$\R$, of the values %
\[%
L^{1,z}\left(\frac{1}{2},1;\st,\sym^m\right)%
\text{ and } 
L^{1,z}\left(1,1;\sym^2,\sym^m\right)%
\]
in the following proposition. Denote by  $\gamma^*$ the constant determined by
\[%
\sum_{p\leq x}\frac{1}{p}=\log_2 x+\gamma^*+O\left(\frac{1}{\log x}\right)\quad (x\geq 2).
\]
If $\gamma$ is the Euler constant, we have
\begin{equation}\label{eq:gamet}
\gamma^*=\gamma+\sum_{p\in\prem}\left[\log\left(1-\frac{1}{p}\right)+\frac{1}{p}\right].
\end{equation}
\begin{propint}\label{propint:asymp}
Let $m\in\{1,2,4\}$.
As $r\to+\infty$ in $\R$, the following estimates hold:
\[%
\log L^{1,\pm r}\left(\frac{1}{2},1;\st,\sym^m\right) =\sym^m_{\pm}r\log_2 r
+\sym^{m,1}_{\pm}r
+O_m\left(\frac{r}{\log r}\right)
\]
and
\[%
\log L^{1,\pm r}\left(1,1;\sym^2,\sym^m\right) = 
\sym^m_{\pm}r\log_2 r
+\sym^{m,1}_{\pm}r
+O_m\left(\frac{r}{\log r}\right)
\]
where
\begin{equation}\label{eq:defsy}
\sym^m_{\pm}:= \max_{g\in\su}\pm\chi_{\sym^m}(g)
\end{equation}
and
\begin{multline}\label{eq:defsyu}
\sym^{m,1}_{\pm}
:= 
\gamma^*\sym^m_{\pm}+\\
\sum_{p\in\prem}
\left\{
\pm\log\left(\pm\max_{g\in\su}\pm D(p^{-1},\sym^m,g)\right)-\frac{\sym^m_{\pm}}{p}
\right\}.
\end{multline}
\end{propint}
\begin{remint}
Some values of $\sym^m_{\pm}$ and $\sym^{m,1}_{\pm}$ may be easily computed %
(see table~\ref{tab_valsym}).
\begin{table}[ht]\label{tab_valsym}
\begin{center}
\begin{tabular}{|c|c|c|c|c|}
\hline
$m$              &  $2$                    & $4$       & even & odd \\
\hline
$\sym^m_{+}$     &  $3$                    & $5$       & $m+1$         & $m+1$ \\
\hline
$\sym^m_{-}$     &  $1$                    & $5/4$     &             & $m+1$\\
\hline
$\sym^{m,1}_{+}$ &  $3\gamma$              & $5\gamma$ & $(m+1)\gamma$ & $(m+1)\gamma$\\
\hline
$\sym^{m,1}_{-}$ &  $\gamma-2\log\zeta(2)$ &         &             & $(m+1)[\gamma-\log\zeta(2)]$ \\ 
\hline
\end{tabular}
\end{center}
\caption{Some values of $\sym^m_{\pm}$ and $\sym^{m,1}_{\pm}$}
\end{table}
The reason why $\sym^{m}_{-}$ is easy computed in the case $m$ odd but not in the 
case $m$ even is that the minimum of the Chebyshev polynomial (see \eqref{eq:trtch}) of second kind is well known
when $m$ is odd (due to symmetry reasons) and not when $m$ is even. 
For $\sym^{m,1}_{-}$, see also Remark~\ref{rem:reason}.
Cogdell \& Michel \cite[Theorem 1.12]{CoMi04} found the same asymptotic behaviour for the non twisted moments.
\end{remint}
Since $L(1/2,f)\geq 0$, we may deduce extremal values of $L(1,\sym^mf)$ with the extra condition of
nonvanishing of $L(1/2,f)$.
\begin{corint}\label{lem:dblcdn}
Let $m\in\{1,2,4\}$ and $N\in\Ncri{\log^{3/2}}$.
Then there exists  $f_m\in\pr$ 
and $g_m\in\pr$
satisfying
\[%
L(1,\sym^{m}f_m)
\geq\eta_+(m)
\left[\log_2(3N)\right]^{\sym^m_{+}}
\qquad\text{and}\qquad
L\left(\frac{1}{2},f_m\right)>0,
\]
\[%
L(1,\sym^{m}g_m)
\leq\eta_-(m)
\left[\log_2(3N)\right]^{-\sym^m_{-}}
\hskip 4,8mm\text{and}\qquad
L\left(\frac{1}{2},g_m\right)>0,
\]
where $\eta_{\pm}(m)=[1+o_m(1)]\exp(\sym^{m,1}_{\pm})$.
\end{corint}

\begin{remint}
The hypothesis $N\in\Ncri{\log^{3/2}}$ is certainly crucial since we can prove
the following result. Fix $m\in\{1,2,4\}$. Denote, for all $\omega\in\Z_{>0}$, by $N_\omega$
the product of the first $\omega$ primes. Assume Grand Riemman hypothesis for the
$m$th symmetric power $L$-functions of primitive forms.
Then, there exist $A_m>0$ and $B_m>0$ such that, for all $\omega\in\Z_{>0}$ and
$f\in\bigcup_{\omega\in\Z_{>0}}\prpn{2}{N_\omega}$ we have
\[%
A_m\leq L(1,\sym^mf)\leq B_m.
\] 
\end{remint}

\subsection{Hecke eigenvalues}\label{sec_hecke}
The Sato-Tate conjecture %
predicts that the sequence of the Hecke eigenvalues at prime numbers %
of a fixed primitive form %
is equidistributed for the Sato-Tate measure %
on $[-2,2]$. 
More precisely, for all $[a,b]\subset[-2,2]$, it is %
expected that %
\begin{equation}\label{eq_ST}%
\lim_{x\to +\infty}\frac{%
\#\left\{%
p\in\prem \colon %
\text{$p\leq x$ and $\lambda_f(p)\in[a,b]$}%
\right\}
}{\#\{p\in\prem \colon p\leq x\}}%
=%
F_{\mathrm{ST}}(b)-F_{\mathrm{ST}}(a)%
\end{equation}
with %
\[%
F_{\mathrm{ST}}(u):= %
\frac{1}{\pi}\int_{-2}^u\sqrt{1-\frac{t^2}{4}}\dd t.
\]
Note that in \eqref{eq_ST}, the primitive form $f$ is fixed %
and hence, the parameter $x$ can not depend on the level of $f$. %
The Sato-Tate conjecture \eqref{eq_ST} is sometimes termed %
\emph{horizontal} Sato-Tate equidistribution conjecture %
in opposition to the \emph{vertical} Sato-Tate equidistribution %
Theorem (due to Sarnak \cite{sar87}, see also \cite{ser97}) in %
which the equidistribution is proved for a fixed prime number $p$. %
For all $[a,b]\subset[-2,2]$, it is proved that%
\[%
\lim_{N\to+\infty}\frac{\left\{%
f\in\pr \colon \lambda_f(p)\in[a,b]%
\right\}}{\#\pr}%
=%
F_{\mathrm{ST}}(b)-F_{\mathrm{ST}}(a).%
\] 
In vertical and horizontal distributions, there %
should be less Hecke eigenvalues in an interval %
near $2$ than in an interval of equal length around %
$0$. %
In Propositions~\ref{propint:rubin} and~\ref{propint:tram}, %
we show that, for many primitive forms, the first few %
(in term of the level) Hecke eigenvalues concentrate %
near (again in term of the level) $2$.
To allow comparisons, we recall the following estimate:%
\[%
\sum_{%
p\leq [\log(2N)]^\epsilon%
} 
\frac{1}{p}
= \log_3(20N) 
\left\{%
1 + O_{\epsilon}\left(\frac{1}{\log_3(20N)}\right)%
\right\}.
\]

Let $N\in\Ncri{\log^{3/2}}$.
For $C>0$, denote by \[\prgde{m}\] the set of primitive forms $f\in\pr$ such that
\begin{equation}\label{eq:dhp}
L(1,\sym^{m}f)
\geq C\left[\log_2(3N)\right]^{\sym^m_{+}}.
\end{equation}
For $C>0$ small enough, such a set is not empty (by an easy adaptation of
\cite[Corollary 1.13]{CoMi04}) and by the method developed in \cite{RW03} its size
is large (although not a positive proportion of $\#\pr$).

\begin{propint}\label{propint:rubin}
Let $m\in\{1,2,4\}$ and $N$ an integer of $\Ncri{\log^{3/2}}$.
For all $\epsilon>0$ and $\xi(N)\to\infty$ $(N\to\infty)$ with
$\xi(N)\leq \log_3(20N)$, for all $f\in\prgde{m}$ 
such that Grand Riemann Hypothesis is true for $L(s,\sym^mf)$, the following estimate holds:
\[%
\sum_{\substack{%
p\leq [\log(2N)]^\epsilon
\\
\lambda_{f}(p^m)\geq\sym^m_+-\xi(N)/\log_3(20N)
}} 
\frac{1}{p}
= \log_3(20N) 
\left\{
1 + O_{\epsilon, m}\left(\frac{1}{\xi(N)}\right)
\right\}.
\]
\end{propint}

Our methods allow to study the small values of the Hecke eigenvalues.
Denote by $\prpte{m}$ the set of primitive forms $f\in\pr$ such that
\[%
L(1,\sym^{m}f)
\leq C\left[\log_2(3N)\right]^{-\sym^m_{-}}.
\]
\begin{propint}\label{propint:tram}
Let $N\in\Ncri{\log^{3/2}}$. For all $\epsilon>0$ and $\xi(N)\to\infty$ $(N\to\infty)$ with
$\xi(N)\leq \log_3(20N)$, for all $f\in\prpte{2}$ 
such that Grand Riemann Hypothesis is true for $L(s,\sym^2f)$, the following estimate holds:
\[%
\sum_{\substack{%
p\leq [\log(2N)]^\epsilon
\\
\lambda_{f}(p)\leq[\xi(N)/\log_3(20N)]^{1/2}
}} 
\frac{1}{p}
= \log_3(20N) 
\left\{
1 + O_{\epsilon}\left(\frac{1}{\xi(N)}\right)
\right\}.
\]
\end{propint}
\begin{remint}
\begin{enumerate}
\item Propositions~\ref{propint:rubin} and~\ref{propint:tram} are also true  with the extra condition $L(1/2,f)>0$.
\item The study of extremal values of symmetric power $L$-functions at $1$
and Hecke eigenvalues in the weight aspect has been done in \cite{LW04} by Lau \& the second author.
\end{enumerate}
\end{remint}

\subsection{Simultaneous extremal values}

Recall that assuming Grand Riemann Hypothesis for $m$th symmetric power $L$-functions, there exists two constants $D_m,D_m'>0$ such that for all
$f\in\pr$, we have
\begin{equation*}
D_m[\log_2(3N)]^{-\sym^m_-}
\leq
L(1,\sym^mf)
\leq
D_m'[\log_2(3N)]^{\sym^m_+}
\end{equation*} 
(see \cite[(1.45)]{CoMi04}). %
We established in Section~\ref{sec_hecke} a link between the %
extremal values of $L(1,\sym^mf)$ and the extremal values %
of $\lambda_f(p^m)$. If we want to study the simultaneous %
extremality of the sequence %
\[%
L(1, \sym^2f),\dotsc, L(1, \sym^{2\ell}f)%
\]
(as $f$ varies), %
we can study the simultaneous %
extremality of the sequence %
\[%
\lambda_f(p^2),\dotsc,\lambda_f(p^{2\ell}).
\]
This is equivalent to the simultaneous %
extremality of the sequence of Chebyshev polynomials %
\[X_2,\dotsc,X_{2\ell}\]
(defined in \eqref{eq:trtch}). But those polynomials are not %
minimal together. An easy resaon is the Clebsh-Gordan %
relation %
\[%
X_\ell^2=\sum_{j=0}^\ell X_{2j}%
\] 
(see \eqref{eq_cg}): the minimal value of the right-hand side would %
be negative if the Chebyshev polynomials were all minimal %
together. Hence, we concentrate on  $L(1,\sym^2f)$ and 
$L(1,\sym^4f)$ and prove that $L(1,\sym^2f)$ and 
$L(1,\sym^4f)$ can not be minimal together but are maximal together.
\begin{propint}\label{propint:simul}
Assume Grand Riemann Hypothesis for symmetric square and fourth symmetric power $L$-functions.
Let $C>0$. 
\begin{enumerate}
\item
There exists no $N\in\Ncri{\log}$ for which there exists $f\in\pr$ satisfying simultaneously
\[%
L(1,\sym^{2}f)
\leq C
\left[\log_2(3N)\right]^{-\sym^2_{-}}
\]
and
\[%
L(1,\sym^{4}f)
\leq C
\left[\log_2(3N)\right]^{-\sym^4_{-}}.
\]
\item
Let $N\in\Ncri{\log}$. 
If $f\in\pr$ satisfies
\[%
L(1,\sym^{2}f)
\geq C
\left[\log_2(3N)\right]^{\sym^2_{+}}
\]
then
\[%
L(1,\sym^{4}f)
\geq C
\left[\log_2(3N)\right]^{\sym^4_{+}}.
\]
\end{enumerate}
\end{propint}

\begin{propint}\label{propint:simulde}
Let $m\geq 1$.
Assume Grand Riemann Hypothesis for symmetric square and $m$th symmetric power $L$-functions.
Let $C,D>0$. 
There exists no $N\in\Ncri{\log}$ for which there exists $f\in\pr$ satisfying simultaneously
\[%
L(1,\sym^{m}f)
\geq C
\left[\log_2(3N)\right]^{\sym^m_{+}}
\]
and
\[%
L(1,\sym^{2}f)
\leq D
\left[\log_2(3N)\right]^{-\sym^2_{-}}.
\]
\end{propint}

\subsection{A combinatorial interpretation of the twisted moments}

The negative moments of $L(1,\sym^2f)$ twisted by $L(1/2,f)$ have a combinatorial interpretation which leads
to Corollary~\ref{corint:identity}. Interpretations of the same flavour have been given in \cite{roy03} and \cite{MR2139690}. 
An interpretation of the traces of Hecke operators, implying the same objects, is also to be found in \cite{MR2114776}.  
We shall denote the vectors with boldface letters:
$\vec{\alpha}=(\alpha_1,\dotsm,\alpha_n)$. Define
$\tr\vec{\alpha}=\sum_{i=1}^n\alpha_i$
and
$\abs{\vec{\alpha}}=\prod_{i=1}^n\alpha_i$.
Let $\mu$ be the Moebius function.
Suppose $n\in\mathbb{N}$ and define 
\[%
\mathcal{E}_n(\vec{b})
:= 
\left\{
\vec{d}\in\N^{n-1}\colon d_i\mid\left(\frac{b_1\dotsm b_i}{d_1\dotsm d_{i-1}},b_{i+1}\right)^2,\,\forall i\in[1,n-1]
\right\},
\]
\[%
w_{-n}(r)
=
\sum_{\substack{\vec{a},\vec{b},\vec{c}\in\N^n\\\abs{\vec{a}\vec{b}^2\vec{c}^3}=r}}
\left[
\prod_{i=1}^n\mu(a_ib_ic_i)\mu(b_i)
\right]
\sum_{\vec{d}\in\mathcal{E}_n(\vec{ab})}\frac{\abs{\vec{d}}}{\abs{\vec{ab}}}
\]
and
\[%
W_{-n}:= 
\prod_{p\in\prem}\sum_{\nu=0}^{+\infty}\frac{w_{-n}(p^\nu)}{p^\nu}.
\]
Using the short expansions of $L(1,\sym^2f)$ (see \eqref{eq:omega}) and $L(1/2,f)$
(see \eqref{eq:repd}) with Iwaniec, Luo \& Sarnak trace formula (see Lemma~\ref{lem:ils}) 
we obtain
\[%
\lim_{\substack{N\to+\infty\\ N\in\Ncri{\log}}}\sum_{f\in\pr}\omegam(f)L\left(\frac{1}{2},f\right)L(1,\sym^2f)^{-n}
=
\zeta(2)^{-n}
W_{-n}.
\]
The method developed in \cite[\S 2.1]{roy03} leads to the following lemma.
\begin{lemint}\label{lem:expcom}
Let $n\geq 0$ and $k\in[0,n]$ be integers. Define
\[%
R_k(p):= 
\begin{cases}
p &\text{if $k=0$ ;}\\
1 &\text{if $k=1$ ;}\\
\displaystyle{\sum_{\substack{\vec{\delta}\in\{-1,0,1\}^{k-1}\\
      \delta_1+\dotsm+\delta_i\leq\max(0,\delta_i)}}p^{\tr\vec{\delta}}} &\text{if $k\geq 2$.}
\end{cases}
\]
Then,
\[%
W_{-n}
=
\frac{1}{\zeta(3)^n}
\prod_{p\in\prem}
\frac{1}{p}
\sum_{k=0}^n
(-1)^k
\binom{n}{k}
R_k(p)
\left(\frac{p}{p^2+p+1}\right)^k.
\]
\end{lemint}
Assume $k\geq 1$.
Writing
\[%
R_k(p)=: \sum_{q=-(k-1)}^1\xi_{k,q}p^q,
\]
the integer $\xi_{k,q}$ is the number of paths in $\Z^2$ which
\begin{itemize}
\item rely $(0,0)$ to $(k-1,q)$ 
\item with steps $(1,-1)$, $(1,0)$ or $(1,1)$ 
\item never going above the abscissas axis
\item \emph{except eventually} with a step $(1,1)$ that is immediately followed
  by a step $(1,-1)$ if it is not the last one.   
\end{itemize}

\begin{figure}[ht!]
\centering
 \includegraphics{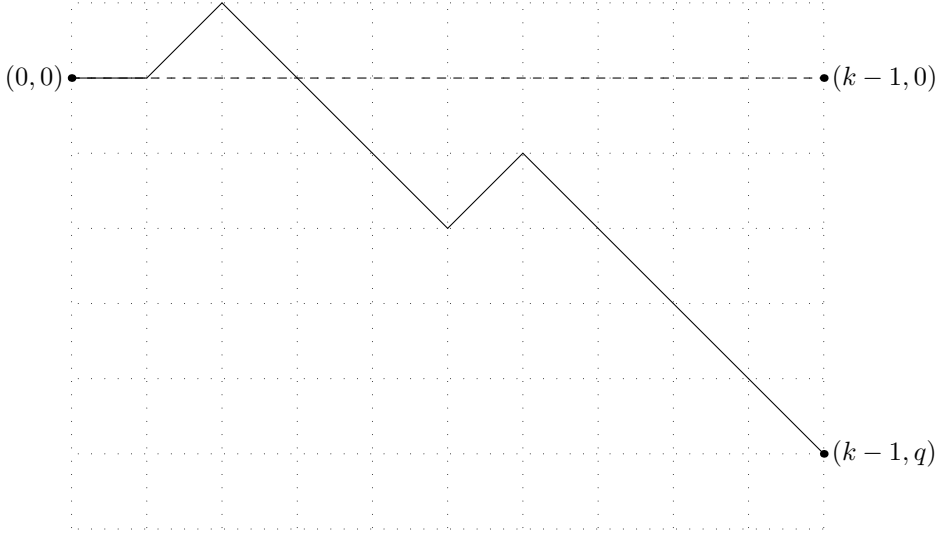}
\caption{A partial Riordan path}\label{fig:rpa}
\end{figure}

In other words, we count \emph{partial Riordan paths} (see
figure~\ref{fig:rpa}). 

\begin{figure}[ht!]
\centering
 \includegraphics{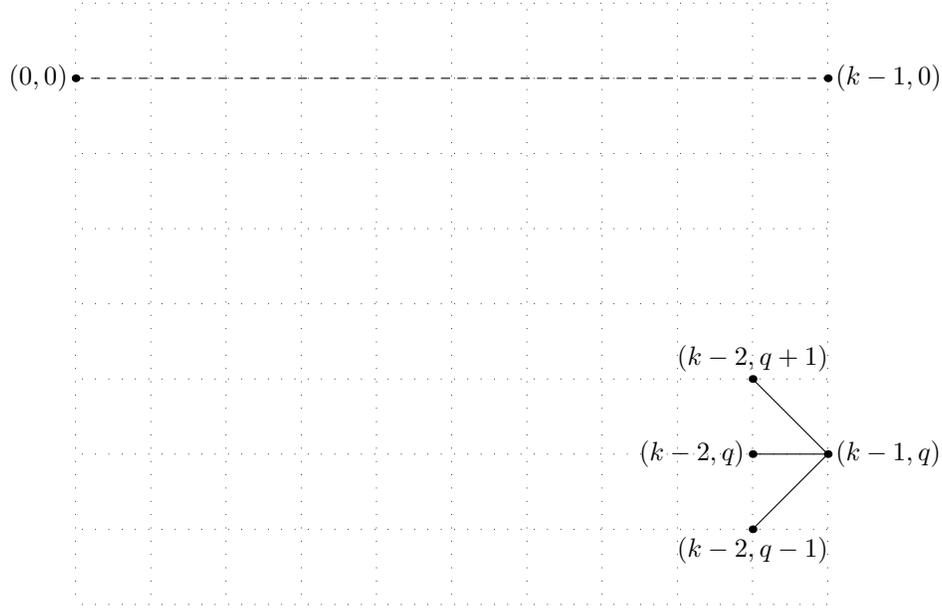}
\caption{Relation between $\xi_{k,q}$, $\xi_{k-1,q-1}$,
  $\xi_{k-1,q}$ and $\xi_{k-1,q+1}$}\label{fig:deux}
\end{figure}

For $q=0$, we obtain a Riordan path. Riordan paths
have been studied in \cite[\S 1.2]{roy03} where the
number of Riordan paths from $(0,0)$ to $(k,0)$ was denoted by $R_{k+2}$ (this number is called the $k+2$th Riordan 
number). We then have
\begin{equation*}
\xi_{k,0}=R_{k+1}.
\end{equation*}
This remains true for $k=0$ since $R_1=0$. The Riordan paths rely to our problem since the first author proved in 
\cite[Proposition 11]{roy03} that
\begin{equation}\label{eq:mmm}
\lim_{\substack{N\to+\infty\\ N\in\Ncri{\log}}}\sum_{f\in\pr}\omegam(f)L(1,\sym^2f)^{-n}
=
\frac{1}{\zeta(3)^n}\prod_{p\in\prem}\ell_n\left(\frac{p}{p^2+p+1}\right)
\end{equation}
where
\begin{align*}
\ell_n(x)&:= 
\sum_{k=0}^n(-1)^k\binom{n}{k}R_kx^k
\\
&=\frac{4}{\pi}\int_0^{\pi/2}\left[1+x(1-4\sin^2\theta)\right]^n\cos^2\theta\dd\theta.
\end{align*}
Using the recursive relation
\[%
R_k(p)
=
\left(p+1+\frac{1}{p}\right)R_{k-1}(p)
-p(p+1)R_{k-1}
\]
which expresses that a path to $(k-1,q)$ has is last step coming from one of the three points 
$(k-2,q+1)$, $(k-2,q)$, $(k-2,q-1)$ (see figure~\ref{fig:deux}) we get
\begin{multline}\label{eq:bruit}
\sum_{k=0}^{n+1}(-1)^k\binom{n+1}{k}R_k(p)\left(\frac{p}{p^2+p+1}\right)^k
\\
=%
\frac{p^2(p+1)}{p^2+p+1}\ell_n\left(\frac{p}{p^2+p+1}\right).
\end{multline}

Reintroducing \eqref{eq:bruit} in Lemma~\ref{lem:expcom} and comparing with \eqref{eq:mmm}
gives
\begin{multline*}
\lim_{\substack{N\to+\infty\\ N\in\Ncri{\log}}}\sum_{f\in\pr}\omegam(f)L\left(\frac{1}{2},f\right)L(1,\sym^2f)^{-n}
=\\
\lim_{\substack{N\to+\infty\\ N\in\Ncri{\log}}}
\frac{1}{\zeta(2)}\sum_{f\in\pr}\omegam(f)L(1,\sym^2f)^{-n+1}.
\end{multline*}
\subsection{A few notation}\label{pa:us}
In this text we shall use the following notation not yet introduced.
We give at the end of the text (see Section~\ref{sec:index}) an index of notation.
If $a$ and $b$ are two complex numbers, then $\delta(a,b)=1$ if $a=b$ and $\delta(a,b)=0$
otherwise. If $n$ is an integer, define $\square(n)=1$ if $n$ is a square and $\square(n)=0$ otherwise.
Remark that $\square$ is \emph{not} the function $\square_1$ (since $\square_1(n)=\delta(n,1)$). 
If $p$ is a prime number,
$v_p(n)$ is the $p$-valuation of $n$. Moreover, if $N$ is
another integer, then we decompose $n$ as
$n=n_Nn^{(N)}$ with $p\mid n_N\Rightarrow p\mid N$ and $(n^{(N)},N)=1$.
The functions $\un_N$ and $\un^{(N)}$ are defined by
\begin{equation}\label{eq:unb}
\un_N(n):= \begin{cases}
1 & \text{if the prime divisors of $n$ divide $N$}\\
0 & \text{otherwise}
\end{cases}
\end{equation}
and
\begin{equation}\label{eq:unh}
\un^{(N)}(n):= \begin{cases}
1 & \text{if $(n,N)=1$}\\
0 & \text{otherwise.}
\end{cases}
\end{equation}
The letters $s$ and $\rho$ are devoted to complex numbers and we set $\Re s=\sigma$ and $\Re\rho=r$.

\section{Modular tools}\label{sec:tools}
In this section, we establish some results needed for the forthcoming proofs of our results.

\subsection{Two standard hypothesis}\label{sec:tsh}
We introduce two standard hypothesis that shall allow us to prove our results for each symmetric power $L$--function.
If $f\in\pr$, we have defined $L(s,\sym^{m}f)$ in \eqref{eq:defsym} as being an Euler product of degree $m+1$.
These representations allow to express the multiplicativity relation of
$n\mapsto\lambda_f(n)$: this function is multiplicative and, if $p\nmid N$ and $\nu\geq 0$, we have
\begin{equation}\label{eq:multhecke}
\lambda_f(p^\nu)=\chi_{\sym^\nu}[g(\theta_{f,p})].
\end{equation}
Recall also that $n\mapsto\lambda_{f}(n)$ is strongly 
multiplicative on integers having their prime factors in the support of $N$ and that if $n\mid N$, then
\begin{equation}\label{eq:krivine}
\abs{\lambda_f(n)}=\frac{1}{\sqrt{n}}.
\end{equation}
The first hypothesis on the automorphy of $L(s,\sym^mf)$
for all $f\in\pr$ is denoted by $\sym^m(N)$. It is has been proved
in the cases $m\in\{1,2,3,4\}$ (see \cite{gej78}, \cite{KS02}, \cite{KS02b} and \cite{kim03}). 
The second hypothesis is concerned with the eventual Landau-Siegel
zero of the $m$th symmetric power $L$-functions, it is denoted by $LSZ^m(N)$ and has been proved
for $m\in\{1,2,4\}$ (see \cite{hol94}, \cite{ghl94}, \cite{hor95} and \cite{raw03}).

Fix $m\geq 1$ and $N$ a squarefree positive integer.
\begin{HypSym}
For every $f\in\pr$,
there exists an automorphic cuspidal selfdual representation of $\gl_{m+1}(\mathbb{A}_{\mathbb{Q}})$
whose local $L$ factors agree with the ones of the function $L(s,\sym^mf)$.
Define
\begin{multline*}
L_\infty(s,\sym^m f)
:= \\
\begin{cases}
\displaystyle{%
\pi^{-s/2}\Gamma\left(\frac{s}{2}\right)2^u\prod_{j=1}^u(2\pi)^{-s-j}\Gamma(s+j)
}
& \text{if $m=2u$ with $u$ even}\\
\displaystyle{%
\pi^{-(s+1)/2}\Gamma\left(\frac{s+1}{2}\right)2^u\prod_{j=1}^u(2\pi)^{-s-j}\Gamma(s+j)
}
& \text{if $m=2u$ with $u$ odd}\\
\displaystyle{%
2^{u+1}
\prod_{j=0}^u(2\pi)^{-s-j-1/2}\Gamma\left(s+j+\frac{1}{2}\right)
}
& \text{if $m=2u+1$.}
\end{cases}
\end{multline*}
Then there exists $\epsilon(\sym^mf)\in\{-1,1\}$ such that
\begin{multline*}
N^{ms/2}
L_\infty(s,\sym^m f)L(s,\sym^m f)
=\\
\epsilon(\sym^mf)
N^{m(1-s)/2}
L_\infty(1-s,\sym^m f)L(1-s,\sym^m f).
\end{multline*}
\end{HypSym}
We refer to \cite{CoMi04} for a discussion on the analytic implications of this conjecture.
The second hypothesis we use is the non existence of Landau-Siegel zero. Let $N$ squarefree 
such that hypothesis $\sym^m(N)$ holds.
\begin{HypLS}
There exists a 
constant $A_m>0$ depending only on $m$ such that for every $f\in\pr$, $L(s,\sym^mf)$
has no zero on the real interval $[1-A_m/\log(2N),1]$.
\end{HypLS}

\subsection{Dirichlet coefficients of the symmetric power $L$-functions}
In this section, we study the Dirichlet coefficients of $L(s,\sym^{m}f)^z$. We derive our study
from the one of Cogdell \& Michel but try to be more explicit in our specific case.
We begin with the polynomial $D$ introduced in \eqref{eq:fato}.
Since $\sym^m$ is selfdual, we have, 
$
D(X,\sym^m,g)\in\R[X]
$
and for $x\in[0,1[$,
\begin{equation}\label{eq:lisyn}
\left(1+{x}\right)^{-m-1}\leq
D(x,\sym^m,g)\leq   \left(1-{x}\right)^{-m-1}.
\end{equation}
\begin{rem}\label{rem:reason}
Note that the upper bound is optimal since the equation %
\[%
\sym^m g=I%
\]
admits always $I$ as a solution whereas the lower bound is optimal only for odd $m$ since
$\sym^m g=-I$ has a solution only for odd $m$.
\end{rem}
Evaluating \eqref{eq:lisyn} at $g=g(\pi)$, we find
\[%
\min_{g\in\su}D\left(X,\sym^{2m+1},g\right)=(1+X)^{-2m-2}.
\]
Next,
\begin{multline*}
D\left[X,\sym^{2m},g\left(\frac{\pi}{2m}\right)\right]
\\=
(1-X)^{-1}\prod_{j=1}^m
\left(
1-Xe^{2j\frac{\pi i}{2m}}
\right)^{-1}
\left(
1-Xe^{-2j\frac{\pi i}{2m}}
\right)^{-1}
\\=
(1+X)^{-1}(1-X^{2m})^{-1}
\end{multline*}
so that
\[%
\min_{g\in\su}D\left(X,\sym^{2m},g\right)\leq (1+X)^{-1}(1-X^{2m})^{-1}.
\]

For every $g\in\su$, define $\lambda_{\sym^m}^{z,\nu}(g)$ by the expansion
\begin{equation}\label{eq:devD}
D(X,\sym^m,g)^z
=: 
\sum_{\nu=0}^{+\infty}\lambda_{\sym^m}^{z,\nu}(g)X^\nu.
\end{equation}
The function $g\mapsto\lambda_{\sym^m}^{z,\nu}(g)$ is central so that it may be expressed as a linear
combination of the characters of irreducible representations of $\su$. These characters are 
defined on the conjugacy classes of $\su$ by
\begin{equation}\label{eq:trtch}
\chi_{\sym^m}[g(\theta)]=\tr\sym^m[g(\theta)]=\frac{\sin[(m+1)\theta]}{\sin\theta}=X_m(2\cos\theta)
\end{equation}
where $X_m$ is the $m$th Chebyshev polynomial of second kind on $[-2,2]$. 
We then have
\begin{equation}\label{eq:supbou}
\lambda_{\sym^m}^{z,\nu}(g)
=
\sum_{m'\geq 0}\mu_{\sym^m,\sym^{m'}}^{z,\nu}\chi_{\sym^{m'}}(g)
\end{equation}
with
\begin{align}
\label{eq:ei}\mu_{\sym^m,\sym^{m'}}^{z,\nu}&=
\int_{\su}
\lambda_{\sym^m}^{z,\nu}(g)
\chi_{\sym^{m'}}(g)\dd g\\
\label{eq:eib} &=
\frac{2}{\pi}\int_0^\pi
\lambda_{\sym^m}^{z,\nu}[g(\theta)]
\sin[(m'+1)\theta]\sin\theta\dd\theta.
\end{align}
We call $\mu_{\sym^m,\sym^{m'}}^{z,\nu}$ the harmonic of $\lambda_{\sym^m}^{z,\nu}$
of order $m'$.
In particular, 
\begin{equation}\label{j25}
\mu_{\sym^m,\sym^{m'}}^{z,0}=\delta(m',0)
\end{equation}
and, since $\lambda_{\sym^m}^{z,1}(g)=z\chi_{\sym^m}(g)$, we have
\begin{equation}\label{eqj26}
\mu_{\sym^m,\sym^{m'}}^{z,1}=z\delta(m,m').
\end{equation}
From the expansion
\begin{equation}\label{eq:paques}
(1-x)^{-z}=\sum_{\nu=0}^{+\infty}\binom{z+\nu-1}{\nu}x^\nu
\end{equation}
we deduce
\begin{equation*}
D[x,\sym^m,g(\theta)]^z
=
\sum_{\nu=0}^{+\infty}
\left\{
\sum_{\substack{\vec{\nu}\in\N^{m+1}\\\tr{\vec\nu}=\nu}}
\left[\prod_{j=0}^m\binom{z+\nu_{j+1}-1}{\nu_{j+1}}\right]
e^{i\ell(m,\vec{\nu})\theta}
\right\}x^\nu
\end{equation*}
with
\begin{equation}\label{eq:ellem}
\ell(m,\vec{\nu}):= m\nu-2\sum_{k=1}^mk\nu_{k+1}
\end{equation}
and gets
\begin{equation}\label{eq:precouperin}
\lambda_{\sym^m}^{z,\nu}[g(\theta)]
=
\sum_{\substack{\vec{\nu}\in\N^{m+1}\\\tr{\vec\nu}=\nu}}
\left[\prod_{j=0}^m\binom{z+\nu_{j+1}-1}{\nu_{j+1}}\right]
e^{i\ell(m,\vec{\nu})\theta}.
\end{equation}
This function is entire in $z$, then assuming that $z$ in real, using that the left hand side is real
in that case, taking the real part in the right hand side and using analytic continuation we have for
all $z$ complex
\begin{equation}\label{eq:couperin}
\lambda_{\sym^m}^{z,\nu}[g(\theta)]
=
\sum_{\substack{\vec{\nu}\in\N^{m+1}\\\tr{\vec{\nu}}=\nu}}
\left[\prod_{j=0}^m\binom{z+\nu_{j+1}-1}{\nu_{j+1}}\right]
\cos\left[\ell(m,\vec{\nu})\theta\right].
\end{equation}
It follows that \eqref{eq:eib} may be rewritten as
\begin{multline*}
\mu_{\sym^m,\sym^{m'}}^{z,\nu}
=
\frac{2}{\pi}
\sum_{\substack{\vec{\nu}\in\N^{m+1}\\\tr{\vec{\nu}}=\nu}}
\left[\prod_{j=0}^m\binom{z+\nu_{j+1}-1}{\nu_{j+1}}\right]
\\
\times
\int_0^{\pi}
\cos\left[\ell(m,\vec{\nu})\theta\right]
\sin[(m'+1)\theta]\sin\theta\dd\theta
\end{multline*}
that is
\begin{equation}\label{eq:leclerc}
\mu_{\sym^m,\sym^{m'}}^{z,\nu}
=
\frac{1}{2}
\sum_{\substack{\vec{\nu}\in\N^{m+1}\\\tr{\vec{\nu}}=\nu}}
\left[\prod_{j=0}^m\binom{z+\nu_{j+1}-1}{\nu_{j+1}}\right]
\Delta(m,m',\vec{\nu})
\end{equation}
with
\begin{equation}\label{eq:valDel}
\Delta(m,m',\vec{\nu})
=
\begin{cases}
2 &\text{if $\ell(m,\vec{\nu})=0$ and $m'=0$}\\ 
1 &\text{if $\ell(m,\vec{\nu})\pm m'=0$ and $m'\neq 0$}\\
-1 &\text{if $\ell(m,\vec{\nu})\pm m'=\mp 2$}\\
0 & \text{otherwise.}
\end{cases}
\end{equation}
In particular, $\mu_{\sym^m,\sym^{m'}}^{z,\nu}=0$ if $m'>m\nu$ thus
\begin{equation}\label{eq:j23}
\lambda_{\sym^m}^{z,\nu}(g)
=
\sum_{m'=0}^{m\nu}\mu_{\sym^m,\sym^{m'}}^{z,\nu}\chi_{\sym^{m'}}(g).
\end{equation}
Equation \eqref{eq:valDel} also immediately gives
\begin{equation}\label{eq:j27}
\mu_{\sym^{2m},\sym^{2m'+1}}^{z,\nu}=0
\end{equation}
and
\begin{equation*}
\mu_{\sym^{2m+1},\sym^{m'}}^{z,\nu}=0\, \text{if $m'$ and $\nu$ have different parity} 
\end{equation*}
for all $m$ and $m'$.

For $m=1$, we have
\begin{equation}\label{eq:diable}
D[X,\st,g(\theta)]
=
\frac{1}{1-2\cos(\theta)X+X^2}
=
\sum_{\nu=0}^{+\infty}X_\nu(2\cos\theta)X^{\nu}
\end{equation}
hence
$\lambda^{1,\nu}_{\st}(g)=\chi_{\sym^{\nu}}(g)$ for all $g\in\su$.
It follows that
\begin{equation}\label{eq:sacregr}
\mu^{1,\nu}_{\st,\sym^{\nu'}}=\delta(\nu,\nu').
\end{equation}

Now, equation \eqref{eq:couperin} implies
\[%
\abs{\lambda_{\sym^m}^{z,\nu}[g(\theta)]}
\leq
\sum_{\substack{\vec{\nu}\in\N^{m+1}\\\tr{\vec{\nu}}=\nu}}
\left[\prod_{j=0}^m\binom{\abs{z}+\nu_{j+1}-1}{\nu_{j+1}}\right]
=
\lambda_{\sym^m}^{\abs{z},\nu}[g(0)]
\]
and
\[%
\sum_{\nu=0}^{+\infty}\lambda_{\sym^m}^{\abs{z},\nu}[g(0)]X^\nu
=
\det[I-X\sym^m\left(g(0)\right)]^{-\abs{z}}
=
(1-X)^{-(m+1)\abs{z}}
\]
so that
\begin{equation}\label{eq:j212}
\abs{\lambda_{\sym^m}^{z,\nu}[g(\theta)]}
\leq
\binom{(m+1)\abs{z}+\nu-1}{\nu}.
\end{equation}
From \eqref{eq:valDel}, remarking that the first case is incompatible with the second and third ones, 
that the two cases in the second case are incompatible and that the two cases of the third case are incompatible,
we deduce that
\[%
\sum_{m'=0}^{m\nu}\abs{\Delta(m,m',\vec{\nu})}\leq 2
\]
and \eqref{eq:leclerc} gives
\begin{equation}\label{eq:mieuxCoMi}
\sum_{m'=0}^{m\nu}\abs{\mu_{\sym^m,\sym^{m'}}^{z,\nu}}
\leq
\binom{(m+1)\abs{z}+\nu-1}{\nu}.
\end{equation}
This is a slight amelioration of Proposition 2.1 of \cite{CoMi04} in the case of $\su$.
It immediately gives
\begin{equation}\label{eq:j211}
\abs{\mu_{\sym^m,\sym^{m'}}^{z,\nu}}
\leq
\binom{(m+1)\abs{z}+\nu-1}{\nu}.
\end{equation}

To conclude this study, define the multiplicative function $n\mapsto\lambda_{\sym^m f}^z(n)$ by the expansion
\begin{equation}\label{eq:seco}
L(s,\sym^mf)^z=: \sum_{n=1}^{+\infty}\lambda_{\sym^m f}^z(n)n^{-s}.
\end{equation}
For easy reference, we collect the results of the previous lines in the
\begin{prop}\label{prop:resume}
Let $N$ be a squarefree integer, $f\in\pr$ ; let $\nu\geq 0$ and $m>0$ be integers and $z$ be a complex number.
Then
\begin{equation*}
\lambda_{\sym^m f}^z(p^\nu)
=
\begin{cases}
\displaystyle{%
\tau_z(p^\nu)
\lambda_f(p^{m\nu})} & \text{if $p\mid N$ }\\
\phantom{\displaystyle{\binom{z+\nu-1}{\nu}\lambda_f(p^{m\nu})}} & \text{}\\
\displaystyle{\sum_{m'=0}^{m\nu}\mu_{\sym^m,\sym^{m'}}^{z,\nu}\lambda_f(p^{m'})} & \text{if $p\nmid N$.}
\end{cases}
\end{equation*}
Moreover,
\begin{equation*}
\abs{\lambda_{\sym^m f}^z(p^\nu)}
\leq
\tau_{(m+1)\abs{z}}(p^\nu)
\end{equation*}
\begin{align*}
\mu^{1,\nu}_{\st,\sym^{\nu'}}&=\delta(\nu,\nu')\\
\mu_{\sym^m,\sym^{m'}}^{z,0}&=\delta(m',0)\\
\mu_{\sym^m,\sym^{m'}}^{z,1}&=z\delta(m,m')\\
\mu_{\sym^{2m},\sym^{2m'+1}}^{z,\nu}&=0\\
\mu_{\sym^{2m+1},\sym^{m'}}^{z,\nu}&=0 
\text{if $m'$ and $\nu$ have different parity,} 
\end{align*}
and
\begin{equation*}
\sum_{m'=0}^{m\nu}\abs{\mu_{\sym^m,\sym^{m'}}^{z,\nu}}
\leq
\binom{(m+1)\abs{z}+\nu-1}{\nu}.
\end{equation*}
\end{prop}

\begin{proof}
We just need to prove the first equation. Assume that $p\mid N$, then
\[%
\sum_{\nu=0}^{\infty}\lambda_{\sym^mf}^z(p^\nu)p^{-\nu s}=[1-\lambda_f(p^\nu)p^{-s}]^{-z}
\]
and the result follows from \eqref{eq:paques} since $n\mapsto\lambda_{f}(n)$ is strongly 
multiplicative on integers having their prime factors in the support of $N$. 
In the case where $p\nmid N$, we have
\[%
\sum_{\nu=0}^{\infty}\lambda_{\sym^mf}^z(p^\nu)p^{-\nu s}=D[p^{-s},\sym^m,g(\theta_{f,p})]^{-z}
\]
so that the results are consequences of 
\[%
\lambda_{\sym^mf}^z(p^\nu)=\lambda^{z,\nu}_{\sym^m}[g(\theta_{f,p})]
\]
and especially of \eqref{eq:j23} and \eqref{eq:multhecke}.
\end{proof}

We shall need the Dirichlet series
\begin{equation}\label{eq:chri}
W_{m,N}^{z,\rho}(s)=\sum_{n=1}^{+\infty}\frac{\varpi_{m,N}^{z,\rho}(n)}{n^s}
\end{equation}
where $\varpi_{m,N}^{z,\rho}$ is the multiplicative function defined by
\begin{equation}\label{eq:chri2}
\varpi_{m,N}^{z,\rho}(p^\nu)=\begin{cases}
0 & \text{if $p\mid N$}\\
\displaystyle{%
\sum_{m'=0}^{m\nu}\frac{\mu_{\sym^m,\sym^{m'}}^{z,\nu}}{p^{\rho m'}}} & \text{otherwise}
\end{cases}
\end{equation}
for all prime number $p$ and $\nu\geq 1$. 
Similarly, define a multiplicative function
$\widetilde{w}_{m,N}^{z,\rho}$ by
\begin{equation}\label{eq:chri3}
\widetilde{w}_{m,N}^{z,\rho}(p^\nu)=\begin{cases}
0 & \text{if $p\mid N$}\\
\displaystyle{%
\sum_{m'=0}^{m\nu}\frac{\abs{\mu_{\sym^m,\sym^{m'}}^{z,\nu}}}{p^{\rho m'}}} & \text{otherwise.}
\end{cases}
\end{equation}
Using equations \eqref{eqj26} and \eqref{eq:j211}, we have
\begin{multline}\label{eq:j410}
\sum_{\nu=0}^{+\infty}\frac{\abs{\widetilde{w}_{m,N}^{z,\rho}(p^\nu)}}{p^{\sigma\nu}}
\leq
\\
\left(1-\frac{1}{p^\sigma}\right)^{-(m+1)\abs{z}}
-\frac{(m+1)\abs{z}}{p^\sigma}
+
\frac{(m+1)\abs{z}}{p^{\sigma+r}}\left(1-\frac{1}{p^\sigma}\right)^{-(m+1)\abs{z}-1}
\end{multline}
so that the series converges for $\Re s>1/2$ and $\Re s+\Re\rho>1$. We actually have an integral representation.
\begin{lem}\label{lem:j46}
Let $s$ and $\rho$ in $\C$ such that $\Re s>1/2$ and $\Re s+\Re\rho>1$.
Let $N$ be squarefree, then
\[%
W_{m,N}^{z,\rho}(s)
=
\prod_{p\nmid N}\int_{\su}D(p^{-s},\sym^m,g)^zD(p^{-\rho},\st,g)\dd g.
\]
Moreover,
\[%
W_{2m,N}^{z,\rho}(s)
=
\frac{1}{\zeta^{(N)}(4\rho)}
\prod_{p\nmid N}\int_{\su}D(p^{-s},\sym^{2m},g)^zD(p^{-2\rho},\sym^2,g)\dd g.
\]
\end{lem}
\begin{rem}\label{rem:plana}
The key point of Corollary~\ref{corint:identity} is the fact that the coefficients appearing in the series expansion of
$D(X,\sym^{2m},g)$ have only even harmonics -- see equations \eqref{eq:j23} and \eqref{eq:j27}. 
This allows to get the second equation in Lemma~\ref{lem:j46}. It does not seem to have an equivalent for $D(X,\sym^{2m+1},g)$.
Actually, we have
\begin{multline*}
W_{2m+1,N}^{z,\rho}(s)
=
\prod_{p\nmid N}\int_{\su}
[1-p^{-4\rho}+p^{-\rho}(1-p^{-2\rho})\chi_{\st}(g)]
\times\\
D(p^{-s},\sym^{2m+1},g)^zD(p^{-2\rho},\sym^2,g)\dd g
\end{multline*}
and the extra term $p^{-\rho}(1-p^{-2\rho})\chi_{\st}(g)$ is the origin of the fail in obtaining Corollary~\ref{corint:identity} for odd powers.
\end{rem}
Before proving Lemma~\ref{lem:j46}, we prove the following one
\begin{lem}\label{lem:cheb}
Let $g\in\su$, $\ell\geq 2$ an integer and $\abs{X}<1$. Then
\[%
\sum_{k=0}^{+\infty}\chi_{\sym^k}(g)X^k=D(X,\st,g)
\] 
and
\[%
\sum_{k=0}^{+\infty}\chi_{\sym^{k\ell}}(g)X^k=[1+\chi_{\sym^{\ell-2}}(g)X]D(X,\st,g^\ell).
\]
In addition,
\[%
\sum_{k=0}^{+\infty}\chi_{\sym^{2k}}(g)X^k=(1-X^2)D(X,\sym^2,g).
\]
\end{lem}
\begin{proof}
Let $g\in\su$. 
Denote by $e^{i\theta}$ and $e^{-i\theta}$ its eigenvalues. The
first point is equation \eqref{eq:diable}.
If $\ell\geq 2$, with $\xi=\exp(2\pi i/\ell)$, $\lambda=e^{i\theta}$ and $x=2\cos\theta$ we have
\[%
\sum_{\nu=0}^{+\infty}X_{\ell\nu}(x)t^{\ell\nu}
=
\frac{1}{\ell}
\sum_{j=0}^{\ell-1}
\frac{1}{(1-\lambda\xi^jt)(1-\overline{\lambda}\xi^jt)}.
\]
On the other hand,
\[%
\sum_{j=0}^{\ell-1}\frac{1}{1-\lambda\xi^jt}
=
\sum_{j=0}^{\ell-1}
\sum_{n=0}^{+\infty}
\lambda^n\xi^{jn}t^n
=
\frac{\ell}{1-\lambda^\ell t^\ell}
\]
so that
\[%
\sum_{\nu=0}^{+\infty}X_{\ell\nu}(x)t^\nu
=
\frac{
1+\dfrac{\lambda^{\ell-1}-\overline{\lambda}^{\ell-1}}{\lambda-\overline{\lambda}}t}
{1-\left(\lambda^\ell+\overline{\lambda}^\ell\right)t+t^2}.
\]
Since 
\[%
\frac{\lambda^{\ell-1}-\overline{\lambda}^{\ell-1}}{\lambda-\overline{\lambda}}=X_{\ell-2}(x)
\]
we obtain the announced result. In the case $\ell=2$, it leads to
\[%
\sum_{k=0}^{+\infty}\chi_{\sym^{2k}}(g)t^k=\frac{1+t}{(1-\lambda^2t)(1-\overline{\lambda}^2t)}=(1-t^2)D(t,\sym^2,g).
\]
\end{proof}
\begin{proof}[\proofname{} of Lemma~\ref{lem:j46}]
It follows from
\[%
\sum_{m'=0}^{m\nu}\frac{\mu_{\sym^m,\sym^{m'}}^{z,\nu}}{p^{\rho m'}}
=
\sum_{m'=0}^{+\infty}\frac{\mu_{\sym^m,\sym^{m'}}^{z,\nu}}{p^{\rho m'}}
\]
and the expression \eqref{eq:ei} that 
\[%
W_{m,N}^{z,\rho}(s)=\prod_{p\nmid N}\int_{\su}\sum_{\nu=0}^{+\infty}\frac{\lambda_{\sym^m}^{z,\nu}(g)}{p^{\nu s}}
\sum_{m'=0}^{+\infty}\frac{\chi_{\sym^{m'}}(g)}{p^{m'\rho}}\dd g.
\]
The first result is then a consequence of Lemma~\ref{lem:cheb}. Next, we deduce from \eqref{eq:j27} that
\[%
W_{2m,N}^{z,\rho}(s)
=
\prod_{p\nmid N}\sum_{\nu=0}^{+\infty}\frac{1}{p^{\nu s}}
\sum_{m'=0}^{+\infty}\frac{\mu_{\sym^{2m},\sym^{2m'}}^{z,\nu}}{p^{2\rho m'}}
\]
and the second result is again a consequence of Lemma~\ref{lem:cheb}.
\end{proof}
We also prove the
\begin{lem}\label{lem:j45}
Let $m\geq 1$. There exists
$c>0$ such that, for all $N$ squarefree, $z\in\C$, 
$\sigma\in \,]1/2, 1]$ and $r\in[1/2,1]$ we have
\[%
\sum_{n\geq 1} 
\frac{\widetilde{w}_{m,N}^{z,\rho}(n)}{n^s}
\leq
\exp\left[
c(z_m+3)\left(
\log_2(z_m+3) + \frac{(z_m+3)^{(1-\sigma)/\sigma}-1}{(1-\sigma)\log(z_m+3)}
\right)
\right]\]
where
\begin{equation}\label{eq:defzm}
z_m:= (m+1)\min\{n\in\N \colon n\geq \abs{z}\}.
\end{equation}
\end{lem}
\begin{proof}
Equation \eqref{eq:j410} gives
\begin{multline*}
\prod_{p^\sigma\leq z_m+3}
\sum_{\nu\geq 0}
\frac{1}{p^{\nu \sigma}} 
\sum_{0\leq \nu'\leq m\nu} 
\frac{\abs{\mu_{\sym^m, \sym^{\nu'}}^{z, \nu}}}{p^{r\nu'}}
\leq \\
\prod_{p^\sigma\leq z_m+3} 
\left(1-\frac{1}{p^\sigma}\right)^{-z_m-1} 
\left(1 + \frac{z_m}{p^{\sigma+1/2}}\right).
\end{multline*}
Using
\[%
\sum_{p\leq y}\frac{1}{p^\sigma}
\leq 
\log_2y + \frac{y^{1-\sigma}-1}{(1-\sigma)\log y}
\]
valid uniformely for $1/2\leq \sigma\leq 1$ and $y\geq e^2$ (see \cite[Lemme 3.2]{TeW03}) we obtain
\begin{multline*}
\prod_{p^\sigma\leq z_m+3}
\sum_{\nu=0}^{+\infty}\frac{\widetilde{w}_{m,N}^{z,r}(p^\nu)}{p^{\nu\sigma}}
\leq\\
\exp
\left[
c(z_m+3)\left(
\log_2(z_m+3)
+ 
\frac{(z_m+3)^{(1-\sigma)/\sigma}-1}{(1-\sigma)\log(z_m+3)}
\right)
\right].
\end{multline*}
For $p^\sigma>z_m+3$, again by \eqref{eq:j410}, we have
\[%
\sum_{\nu\geq 0}
\frac{1}{p^{\nu \sigma}} 
\sum_{0\leq \nu'\leq m\nu} 
\frac{\abs{\mu_{\sym^m, \sym^{\nu'}}^{z, \nu}}}{p^{r\nu'}}
\leq 
1 + \frac{c(z_m+3)^2}{p^{2\sigma}} 
+ 
\frac{c(z_m+3)}{p^{\sigma+1/2}},\]
so that
\begin{align*}
\prod_{p^\sigma>z_m+3}
\sum_{\nu=0}^{+\infty}\frac{\widetilde{w}_{m,N}^{z,r}(p^\nu)}{p^{\nu\sigma}}
&\leq 
e^{c(z_m+3)^{1/\sigma}/\log(z_m+3)}\\
&\leq
\exp\left[
c(z_m+3)
\frac{(z_m+3)^{(1-\sigma)/\sigma}-1}{(1-\sigma)\log(z_m+3)}
\right].
\end{align*}
\end{proof}
For the primes dividing the level, we have the
\begin{lem}\label{lem:edf}
Let $\ell,m\geq 1$. For $\sigma\in \,]1/2, 1]$ and $r\in[1/2,1]$ we have
\begin{equation*}
\prod_{p\mid N}\int_{\su}D(p^{-s},\sym^m,g)^zD(p^{-\rho},\sym^{\ell},g)\dd g
=\\
1
+
O_{m,\ell}(\Err)
\end{equation*}
with
\[%
\Err:= 
\frac{\omega(N)}{P^-(N)^{2r}}
+
\frac{\abs{z}\omega(N)}{P^-(N)^{r+\sigma}}
+
\frac{\abs{z}^2\omega(N)}{P^-(N)^{2\sigma}}
\]
uniformely for
\begin{equation*}
\begin{cases}
N\in\Ncri{%
\max\bigl\{
\omega(\cdot)^{1/(2r)},\, 
[\abs{z}\omega(\cdot)]^{1/(r+\sigma)},\, 
[\abs{z}^2\omega(\cdot)]^{1/(2\sigma)}
\bigl\}
}
,\text{}\\
z\in\C.
\end{cases}
\end{equation*}
\end{lem}
\begin{proof}
Write
\[%
\Psi^z_{m,\ell}(p):= \int_{\su}D(p^{-s},\sym^m,g)^zD(p^{-\rho},\sym^{\ell},g)\dd g.
\]
Using \eqref{eq:supbou} and the orthogonality of characters, we have
\[%
\Psi^z_{m,\ell}(p)=\sum_{\nu_1=0}^{+\infty}\sum_{\nu_2=0}^{+\infty}p^{-\nu_1s-\nu_2\rho}
\sum_{\nu=0}^{\min(m\nu_1,\ell\nu_2)}\mu^{z,\nu_1}_{\sym^m,\sym^\nu}
\mu^{1,\nu_2}_{\sym^{\ell},\sym^\nu}.
\]
Proposition~\ref{prop:resume} gives
\begin{align*}
\abs{\Psi^z_{m,\ell}(p)-1} &\leq 
\sum_{\nu_2=2}^{+\infty}\binom{\nu_2+\ell}{\nu_2}\frac{1}{p^{r\nu_2}}
+\frac{\abs{z}}{p^{\sigma}}\sum_{\nu_2=1}^{+\infty}\binom{\nu_2+\ell}{\nu_2}\frac{1}{p^{r\nu_2}}\\
&\phantom{\leq\leq}
+
\sum_{\nu_1=2}^{+\infty}\binom{(m+1)\abs{z}+\nu_1-1}{\nu_1}\frac{1}{p^{\sigma\nu_1}}
\sum_{\nu_2=0}^{+\infty}\binom{\nu_2+\ell}{\nu_2}\frac{1}{p^{r\nu_2}}\\
&\ll_{m,\ell}
\frac{1}{p^{2r}}+\frac{\abs{z}}{p^{r+\sigma}}+\frac{\abs{z}^2}{p^{2\sigma}}
\end{align*}
which leads to the result.
\end{proof}
Using \eqref{eq:sacregr} we similarly can prove the
\begin{lem}\label{lem:diable}
Let $m\geq 1$ and $z\in\C$, then
\[%
\int_{\su}D(p^{-1},\sym^m,g)^zD(p^{-1/2},\st,g)\dd g
=
1+O_{m}\left(\frac{\abs{z}}{p^{1+m/2}}\right)
\]
for $p\geq (m+1)\abs{z}+3$.
\end{lem}

\subsection{Dirichlet coefficients of a product of $L$-functions}

The aim of this section is to study the Dirichlet coefficients of the product
\[%
L(s,\sym^2f)L(s,\sym^mf)^z.
\]

Define $\lambda_{\sym^2,\sym^m}^{1,z,\nu}(g)$
for every $g\in\su$ by the expansion
\begin{equation}\label{eq:dvu}
D(x,\sym^2,g)D(x,\sym^m,g)^z=: \sum_{\nu=0}^{+\infty}\lambda_{\sym^2,\sym^m}^{1,z,\nu}(g)x^\nu.
\end{equation}
We have
\begin{equation}\label{eq:noedort}
\lambda_{\sym^2,\sym^m}^{1,z,\nu}(g)=\sum_{\substack{(\nu_1,\nu_2)\in\N^2\\\nu_1+\nu_2=\nu}}
\lambda_{\sym^2}^{1,\nu_1}(g)\lambda_{\sym^m}^{z,\nu_2}(g)
\end{equation}
from that we deduce, using \eqref{eq:j212}, that
\begin{equation*}
\abs{%
\lambda_{\sym^2,\sym^m}^{1,z,\nu}(g)
}
\leq
\binom{(m+1)\abs{z}+2+\nu}{\nu}.
\end{equation*} 
Since $\lambda_{\sym^2,\sym^m}^{1,z,\nu}$ is central, there exists $(\mu^{1,z,\nu}_{\sym^2,\sym^m,\sym^{m'}})_{m'\in\N}$
such that, for all $g\in\su$ we have
\begin{equation}\label{eq:nsle}
\lambda_{\sym^2,\sym^m}^{1,z,\nu}(g)
=
\sum_{m'=0}^{+\infty}\mu^{1,z,\nu}_{\sym^2,\sym^m,\sym^{m'}}\chi_{\sym^{m'}}(g)
\end{equation}
where
\begin{align}
\notag
\mu^{1,z,\nu}_{\sym^2,\sym^m,\sym^{m'}}
&=
\int_{\su}\lambda_{\sym^2,\sym^m}^{1,z,\nu}(g)\chi_{\sym^{m'}}(g)\dd g\\
\label{eq:noepleure}
&=
\frac{2}{\pi}\int_{0}^{\pi}\lambda_{\sym^2,\sym^m}^{1,z,\nu}(g)\sin[(m'+1)\theta]\sin\theta\dd\theta.
\end{align}
The Clebsh-Gordan relation \cite[\S III.8]{vil68} is
\begin{equation}\label{eq_cg}%
\chi_{\sym^{m'_1}}\chi_{\sym^{m'_2}}%
=%
\sum_{r=0}^{\min(m'_1,m'_2)}\chi_{\sym^{m'_1+m'_2-2r}}.%
\end{equation}
In addition with \eqref{eq:noedort} and \eqref{eq:j23}, this relation leads to
\begin{multline}\label{eq:help}
\mu^{1,z,\nu}_{\sym^2,\sym^m,\sym^{m'}}
\\
=%
\sum_{\substack{(\nu_1,\nu_2)\in\N^2\\\nu_1+\nu_2=\nu}}
\mathop{%
\sum_{m'_1=0}^{2\nu_1}
\sum_{m'_2=0}^{m\nu_2}
}_{\substack{\abs{m'_2-m'_1}\leq m'\leq m'_1+m'_2\\ m'_1+m'_2\equiv m'\pmod{2}}}
\mu^{1,\nu_1}_{\sym^2,\sym^{m'_1}}
\mu^{z,\nu_2}_{\sym^m,\sym^{m'_2}}.
\end{multline}
It follows immediately from \eqref{eq:help} that
\[%
\mu^{1,z,\nu}_{\sym^2,\sym^m,\sym^{m'}}=0 \quad \text{if $m'>\max(2,m)\nu$.}
\]
Using also \eqref{j25}, we obtain
\[%
\mu^{1,z,0}_{\sym^2,\sym^m,\sym^{m'}}=\delta(m',0)
\]
and \eqref{eqj26} gives
\[%
\mu^{1,z,1}_{\sym^2,\sym^m,\sym^{m'}}=z\delta(m',m)+\delta(m',2).
\]
Finally, equation\eqref{eq:help} and \eqref{eq:j27} give
\[%
\mu^{1,z,\nu}_{\sym^2,\sym^{2m},\sym^{2m'+1}}=0.
\] 
By equations \eqref{eq:noedort} and \eqref{eq:precouperin} we get
\begin{multline*}
\lambda_{\sym^2,\sym^m}^{1,z,\nu}[g(\theta)]
=\\
\sum_{\substack{(\nu',\nu'')\in\N^2\\\nu'+\nu''=\nu}}
\sum_{\substack{(\vec{\nu'},\vec{\nu''})\in\N^3\times\N^{m+1}\\
\tr\vec{\nu'}=\nu'\\
\tr\vec{\nu''}=\nu''
}}\left[
\prod_{j=0}^m\binom{z+\nu''_{j+1}-1}{\nu''_{j+1}}
\right]
\cos[\ell(2,m;\vec{\nu'},\vec{\nu''})\theta]
\end{multline*}
with
\begin{equation}\label{eq:miap}
\ell(2,m;\vec{\nu'},\vec{\nu''})
=
2\nu'+m\nu''-2\sum_{k=1}^2k\nu'_{k+1}-2\sum_{k=1}^mk\nu''_{k+1}.
\end{equation}
We deduce then from \eqref{eq:noepleure} that
\begin{multline*}
\mu^{1,z,\nu}_{\sym^2,\sym^m,\sym^{m'}}
=\\
\frac{1}{2}
\sum_{\substack{(\nu',\nu'')\in\N^2\\\nu'+\nu''=\nu}}
\sum_{\substack{(\vec{\nu'},\vec{\nu''})\in\N^3\times\N^{m+1}\\
\tr\vec{\nu'}=\nu'\\
\tr\vec{\nu''}=\nu''
}}\left[
\prod_{j=0}^m\binom{z+\nu''_{j+1}-1}{\nu''_{j+1}}
\right]
\Delta(2,m,m';\vec{\nu'},\vec{\nu''})
\end{multline*}
with
\begin{multline}\label{eq:monpey}
\Delta(2,m,m';\vec{\nu'},\vec{\nu''})
\\
:=%
\frac{4}{\pi}\int_0^{\pi}
\cos[\ell(2,m;\vec{\nu'},\vec{\nu''})\theta]
\sin[(m'+1)\theta]\sin\theta\dd\theta.
\end{multline}
From
\[%
\sum_{m'=0}^{\max(2,m)\nu}\abs{\Delta(2,m,m';\vec{\nu'},\vec{\nu''})}\leq 2
\]
we then have
\begin{multline*}
\sum_{m'=0}^{\max(2,m)\nu}%
\abs{\mu^{1,z,1}_{\sym^2,\sym^m,\sym^{m'}}}
\\
\leq%
\sum_{\substack{(\nu',\nu'')\in\N^2\\\nu'+\nu''=\nu}}%
\binom{2+\nu'}{\nu'}\binom{(m+1)\abs{z}+\nu''-1}{\nu''}%
\\
\leq%
\binom{(m+1)\abs{z}+2+\nu''}{\nu''}.
\end{multline*}
To conclude this study, define the multiplicative function \[n\mapsto\lambda_{\sym^2f,\sym^m f}^{1,z}(n)\] by the 
expansion
\begin{equation}\label{eq:conf}
L(s,\sym^2f)L(s,\sym^mf)^z=: \sum_{n=1}^{+\infty}\lambda_{\sym^2f,\sym^m f}^{1,z}(n)n^{-s}.
\end{equation}
The preceding results imply the
\begin{prop}\label{prop:resume2}
Let $N$ be a squarefree integer, $f\in\pr$ ; let $\nu\geq 0$ and $m>0$ be integers and $z$ be a complex number.
Then
\begin{equation*}
\lambda_{\sym^2f,\sym^m f}^{1,z}(p^\nu)
=
\begin{cases}
\displaystyle{%
\sum_{\nu'=0}^{\nu}
\tau_z(p^{\nu'})
\lambda_f(p^{m\nu'})}p^{\nu'-\nu} & \text{if $p\mid N$ }\\
\phantom{\displaystyle{\binom{z+\nu-1}{\nu}\lambda_f(p^{m\nu})}} & \text{}\\
\displaystyle{\sum_{m'=0}^{\max(2,m)\nu}
\mu_{\sym^2,\sym^m,\sym^{m'}}^{1,z,\nu}\lambda_f(p^{m'})} & \text{if $p\nmid N$.}
\end{cases}
\end{equation*}
Moreover,
\begin{equation*}
\abs{\lambda_{\sym^2f,\sym^m f}^{1,z}(p^\nu)}
\leq
\tau_{(m+1)\abs{z}+3}(p^\nu)
\end{equation*}
\begin{align*}
\mu^{1,z,0}_{\sym^2,\sym^m,\sym^{m'}}&=\delta(m',0)\\
\mu^{1,z,1}_{\sym^2,\sym^m,\sym^{m'}}&=z\delta(m',m)+\delta(m',2)\\
\mu^{1,z,\nu}_{\sym^2,\sym^{2m},\sym^{2m'+1}}&=0,
\end{align*}
and
\begin{equation*}
\sum_{m'=0}^{\max(2,m)\nu}\abs{\mu^{1,z,\nu}_{\sym^2,\sym^{m},\sym^{m'}}}
\leq \binom{(m+1)\abs{z}+2+\nu}{\nu}.
\end{equation*}
\end{prop}

\subsection{Trace formulas}%
In this section, we establish a few mean value results for Dirichlet coefficients of the different $L$--functions
we shall encounter.

Let $f\in\pr$.
Denote by $\epsilon_f(N):= \epsilon(\sym^1f)$ the sign of the functional equation satisfied by $L(s,f)$. We have
\begin{equation}\label{eq:tirage}
\epsilon_f(N)=-\mu(N)\sqrt{N}\lambda_f(N)\in\{-1,1\}.
\end{equation}

The following trace formula is due to Iwaniec, Luo \& Sarnak \cite[Corollary 2.10]{ils00}.
\begin{lem}\label{lem:ils}
Let $N\geq 1$ be a squarefree integer
and $m\geq 1$, $n\geq 1$ two integers satisfying $(m,N)=1$ and $(n,N^2)\mid N$.
Then
\begin{equation*}
\sum_{f\in\pr}\omegam(f)
\lambda_f(m)\lambda_f(n)=
\delta(m,n)
\\+
O(\Err)
\end{equation*}
with
\[%
\Err:= 
\frac{\tau(N)^2\log_2(3N)}{N}\,\frac{(mn)^{1/4}\,\tau_3[(m,n)]}{\sqrt{(n,N)}}
\, \log(2mnN).
\]
\end{lem}
We shall need a slightly different version of this trace formula (we actually only remove the condition 
$(n,N)=1$ from \cite[Proposition 2.9]{ils00}).
\begin{lem}\label{lem:aecrire}
Let $N\geq 1$ be a squarefree integer
and $m\geq 1$, $n\geq 1$ two integers satisfying $(m,N)=1$ and $(n,N^2)\mid N$.
Then
\begin{equation*}
\sum_{f\in\pr}\omegam(f)
\left[1+\epsilon_f(N)\right]
\lambda_f(m)\lambda_f(n)=
\delta(m,n)
+
O(\Err)
\end{equation*}
with
\begin{multline*}
\Err:= 
\frac{\delta(n,mN)}{\sqrt{N}}
\\+
\frac{\tau(N)^2\log_2(3N)}{N^{3/4}}
\frac{(mn)^{1/4}}{\sqrt{(n,N)}}\log(2mnN)
\left[
\frac{\tau_3\left[(m,n)\right]}{N^{1/4}}
+
\frac{\tau\left[(m,n)\right]}{\sqrt{(n,N)}}
\right].
\end{multline*}
\end{lem}
\begin{proof}
By Lemma~\ref{lem:ils}, it suffices to prove that
\begin{multline*}
\sum_{f\in\pr}\omegam(f)
\epsilon_f(N)\lambda_f(m)\lambda_f(n)
\ll
\\
\frac{\delta(n,mN)}{\sqrt{N}}
+
\frac{\tau(N)^2\log_2(3N)}{N^{3/4}}
\frac{(mn)^{1/4}}{(n,N)}\tau\left[(m,n)\right]\log(2mnN).
\end{multline*}
Since $\epsilon_f(N)=-\mu(N)\sqrt{N}\lambda_f(N)$, we shall estimate
\[%
R:= \sqrt{N}\sum_{f\in\pr}\omegam(f)
\lambda_f(m)\lambda_f(n)\lambda_f(N).
\]
The multiplicativity relation 
\eqref{eq:multhecke} and equation \eqref{eq:krivine} give 
\begin{align*}
R
&=\sqrt{N}\sum_{f\in\pr}\omegam(f)
\lambda_f(m)\lambda_f(n^{(N)})\lambda_f\left(n_N\right)^2\lambda_f\left(\frac{N}{n_N}\right)\\
&=
\frac{\sqrt{N}}{n_N}
\sum_{d\mid(m,n^{(N)})}
\sum_{f\in\pr}\omegam(f)
\lambda_f\left(\frac{mn^{(N)}}{d^2}\right)
\lambda_f\left(\frac{N}{n_N}\right).
\end{align*}
Then, Lemma~\ref{lem:ils} leads to the result since $MN^{(N)}/d^2=N/n_N$ implies
$N=n_N$, $m=n^{(N)}$ and $d=m$.
\end{proof}

We also prove a trace formula implying the Dirichlet coefficients of the symmetric power $L$-functions.
\begin{lem}\label{lem:fm}
Let $N$ be a squarefree integer,
$(m, n, q)$ be nonnegative integers and $z$ be a complex number.
Then
\[%
\sum_{f\in\pr}\omegam(f)\left[1+\epsilon_f(N)\right] 
\lambda_{\sym^{m}f}^z(n)\lambda_f(q)
= 
w_{m}^{z}(n, q) 
+ 
O(\Err)
\]
with
\begin{equation}\label{eq:defterpr}
w_{m}^{z}(n, q)
:=  
\tau_z(n_N) 
\frac{\square(n_N^mq_N)}{\sqrt{n_N^mq_N}}  
\mathop{%
\prod_{1\leq j\leq r} 
\sum_{0\leq \nu'_j\leq m\nu_{j}}}_{p_1^{\nu'_1}\cdots p_r^{\nu'_r} = q^{(N)}
}
\mu_{\sym^m, \sym^{\nu'_j}}^{z, \nu_{j}}
\end{equation}
where
\[%
n^{(N)}=\prod_{j=1}^rp_j^{\nu_j},\quad (p_1<\dotsm<p_j)
\]
and
\[%
\Err:= 
\frac{\tau(N)^2\log_2(3N)}{N^{3/4}}n^{m/4} 
\tau_{(m+1)\abs{z}}(n)
 \tau(q) q^{1/4} \log(2Nnq).
\]
The implicit constant is absolute.
\end{lem}
\begin{proof}
Let
\[%
S:= \sum_{f\in\pr}\omegam(f)\left[1+\epsilon_f(N)\right] 
\lambda_{\sym^{m}f}^z(n)\lambda_f(q).
\]
Writing $n_N^QM_N=g^2h$ with $h$ squarefree, equation  \eqref{eq:krivine} and
Proposition~\ref{prop:resume} give
\begin{multline*}
S=
\frac{\tau_z(n_N)}{g}
\sum_{(\nu'_i)_{1\leq i\leq r}\in
\bigtimes{i=1}{r}
[0,m\nu_i]}
\left(\prod_{j=1}^{r}\mu_{\sym^m,\sym^{\nu'_j}}^{z,\nu_j}\right)\\
\times
\sum_{d\mid\left(q^{(N)},\prod_{j=1}^{r}p_j^{\nu'_j}\right)}
\sum_{f\in\pr}\omegam(f)[1+\epsilon_f(N)]\lambda_f(h)\lambda_f
\left(\frac{q^{(N)}}{d^2}
\prod_{j=1}^{r}p_j^{\nu'_j}
\right).
\end{multline*}
Then, since $h\mid N$, Lemma~\ref{lem:aecrire} gives $S=P+E$ with
\[%
P=
\frac{\tau_z(n_N)}{g}
\prod_{j=1}^{r} 
\sum_{\nu'_j=0}^{m\nu_j} 
\mu_{\sym^m,\sym^{\nu'_j}}^{z, \nu_j}
\sum_{%
\substack{%
d\mid\left(q^{(N)},\prod_{j=1}^{r}p_j^{\nu'_j}\right)
\\
q^{(N)} p_1^{\nu'_1}\cdots p_r^{\nu'_r}/d^2 = h
}} 1
\]
and
\begin{multline*}
E\ll 
\\
\frac{\tau(N)^2\log_2(3N)}{N^{3/4}}
\frac{n^{m/4}\tau_{\abs{z}}(n_N)}{n_N^{m/2}} 
\frac{q^{1/4}\tau(q)}{q_N^{1/2}}
\frac{\log(2Nnq)}{g^{1/2}}
\prod_{j=1}^{r}  
\sum_{\nu'_j=0}^{m\nu_j}
\abs{\mu_{\sym^m,\sym^{\nu'_j}}^{z, \nu_j}}.
\end{multline*}
Using \eqref{eq:mieuxCoMi}, we obtain
\begin{equation*}
E
\ll
\frac{\tau(N)^2\log_2(3N)}{N^{3/4}}
n^{m/4}
q^{1/4}\tau(q)
\log(2Nnq)
\tau_{(m+1)\abs{z}}(n).
\end{equation*}
We transform $P$ as the announced principal term since $q^{(N)} p_1^{\nu'_1}\cdots p_r^{\nu'_r}/d^2 = h$
implies $p_1^{\nu'_1}\cdots p_r^{\nu'_r}=q^{(N)}=d$ and $h=1$.
\end{proof}

Similarly to Lemma~\ref{lem:fm}, we prove the
\begin{lem}\label{lem:siegfried}
Let $k$, $N$, $m$, $n$ be positive integers, $k$ even, $N$ squarefree. Let $z\in\C$. 
Then
\[%
\sum_{f\in\pr}\omegam(f)\lambda^{1,z}_{\sym^2f,\sym^mf}(n)
=
w^{1,z}_{2,m}(n)
+
O_{k,m}(\Err)
\]
with
\[%
\Err:= 
\frac{\tau(N)^2\log_2(3N)}{N}n^{\max(2,m)\nu/4}r^{1,z}_{2,m}(n)\log(2nN)
\]
where $w^{1,z}_{2,m}$ and $r^{1,z}_{2,m}$ are the multiplicative functions defined by
\[%
w^{1,z}_{2,m}(p^{\nu}):= \begin{cases}
\displaystyle{%
\sum_{\nu'=0}^{\nu}\frac{\tau_{z}(p^{\nu'})\square(p^{m\nu'})}{p^{\nu-\nu'+m\nu'/2}}
} & \text{if $p\mid N$}\\
\phantom{\displaystyle{%
\sum_{\nu'=0}^{\nu}\frac{\tau_{z}(p^{\nu'})\square(p^{m\nu'})}{p^{\nu-\nu'+m\nu'/2}}
}} & \text{}\\
\mu^{1,z,\nu}_{\sym^2,\sym^m,\sym^0} & \text{if $p\nmid N$}
\end{cases}
\]
and
\[%
r^{1,z}_{2,m}(p^{\nu}):= \begin{cases}
\displaystyle{%
\sum_{\nu'=0}^{\nu}\frac{\tau_{\abs{z}}(p^{\nu'})}{p^{\nu-\nu'+m\nu'/2}}
} & \text{if $p\mid N$}\\
\phantom{\displaystyle{%
\sum_{\nu'=0}^{\nu}\frac{\tau_{z}(p^{\nu'})\square(p^{m\nu'})}{p^{\nu-\nu'+m\nu'/2}}
}} & \text{}\\
\binom{(m+1)\abs{z}+\nu+2}{\nu}
& \text{if $p\nmid N$.}
\end{cases}
\]
\end{lem}

\subsection{Mean value formula for the central value of $L(s,f)$}%
Using the functional equation of $L(s,f)$ (see hypothesis $\sym^1(N)$, which is proved in this case)
and contour integrations (see \cite[Theorem 5.3]{IwKo04} for a beautiful explanation) we write
\begin{equation}\label{eq:repd}
L\left(\frac{1}{2},f\right)=\left[1+\epsilon_f(N)\right]\sum_{q=1}^{+\infty}\frac{\lambda_f(q)}{\sqrt{q}}
\exp\left(-\frac{2\pi q}{\sqrt{N}}\right).
\end{equation}
From \eqref{eq:repd} and Lemma~\ref{lem:aecrire} we classically deduce the
\begin{lem}\label{lem:moyLdem}
Let $N$ be a squarefree integer, then
\[%
\sum_{f\in\pr}\omegam(f)L\left(\frac{1}{2},f\right)
=
\zeta_N(2)+O\left(\frac{\tau(N)^2\log(2N)\log_2(3N)}{N^{3/8}}\right).
\]
\end{lem}
\begin{rem}
For $N$ squarefree, we have
\[%
\zeta_N(2)=1+O\left(\frac{\tau(N)}{P^-(N)^2}\right).
\]
\end{rem}
Note that the ``big O'' term may be not small: for all $\omega\geq 1$, let $N_\omega$ be the product of
the $\omega$ first prime numbers, then Mertens Theorem implies that
\[%
\zeta_{N_\omega}(2)\sim\zeta(2)
\]
as $\omega$ tends to infinity. 
\begin{proof}[\proofname{} of Lemma~\ref{lem:moyLdem}]
Equation \eqref{eq:repd} leads to
\begin{multline*}
\sum_{f\in\pr}\omegam(f)L\left(\frac{1}{2},f\right)
=
\\
\sum_{q=1}^{+\infty}\frac{1}{\sqrt{q}}\exp\left(-\frac{2\pi q}{\sqrt{N}}\right)
\sum_{f\in\pr}\omegam(f)
\left[1+\epsilon_f(N)\right]
\lambda_f(q).
\end{multline*}
Writing $q=m\ell^2n$ with $(m,N)=1$, $\ell^2n$ having same prime factors as $N$ and $n$ squarefree, we deduce from
the multiplicativity of $n\mapsto\lambda_f(n)$, its strong multiplicativity of numbers with support included in that
of $N$ and \eqref{eq:krivine} that
\[%
\lambda_f(q)=\frac{1}{\ell}\lambda_f(m)\lambda_f(n).
\]
Then Lemma~\ref{lem:aecrire} gives
\[%
\sum_{f\in\pr}\omegam(f)L\left(\frac{1}{2},f\right)
=
P(N)+O\left(E_1+\tau(N)^2\log_2(3N)(E_2+E_3)\right)
\]
where
\[%
P(N)=\sum_{\ell=1}^{+\infty}\frac{\un_N(\ell)}{\ell^2}\exp\left(-\frac{2\pi\ell^2}{\sqrt{N}}\right)
\]
and
\[%
E_1=\frac{1}{N}\sum_{\ell=1}^{+\infty}\frac{1}{\ell^2}\exp(-2\pi\ell^2\sqrt{N})
\ll
\frac{1}{N},
\]
\begin{align*}
E_2&=
\frac{1}{N}
\sum_{\substack{q=1\\ q=m\ell^2n}}^{+\infty}
\frac{\un_N(\ell n)\un^{(N)}(m)\mu(n)^2\log(2mnN)}{m^{1/4}\ell^2n^{3/4}}\exp\left(-\frac{2\pi m\ell^2n}{\sqrt{N}}\right)\\
&\ll
\frac{1}{N}
\sum_{q=1}^{+\infty}\frac{\log(2qN)}{q^{1/4}}\exp\left(-\frac{2\pi q}{\sqrt{N}}\right)\\
&\ll
\frac{\log(2N)}{N^{5/8}}
\end{align*}
and
\begin{align*}
E_3&=
\frac{1}{N^{3/4}}
\sum_{\substack{q=1\\ q=m\ell^2n}}^{+\infty}
\frac{\un_N(\ell n)\un^{(N)}(m)\mu(n)^2\log(2mnN)}{m^{1/4}\ell^2n^{5/4}}\exp\left(-\frac{2\pi m\ell^2n}{\sqrt{N}}\right)\\
&\ll
\frac{\log(2N)}{N^{3/8}}.
\end{align*}
We conclude by expressing $P(N)$ \textit{via} the inverse Mellin transform of $\exp$ and doing a contour integration
obtaining
\[%
P(N)=\zeta_N(2)+O_\epsilon(N^{-1/2+\epsilon})
\] 
for all $\epsilon>0$.
\end{proof}
 \section{Twisting by $L(1/2,f)$}
The goal of this section is the proof of Theorem~\ref{thmint:A} and Proposition~\ref{propint:B}.
\subsection{\proofname{} of Theorem~\ref{thmint:A}}

Let $z\in\C$ and $x\geq 1$, define
\begin{equation}\label{eq:omega}
\omega_{\sym^mf}^z(x):= \sum_{n=1}^{+\infty}\frac{\lambda_{\sym^m}^z(n)}{n}e^{-n/x}
\end{equation}
for all $f\in\pr$ and proves the
\begin{lem}\label{lem:0918}
Let $N$ be a squarefree integer, $m\in\Z_{>0}$, $x\geq 1$ and $z\in\C$. 
Then
\begin{multline*}
\sum_{f\in\pr}\omegam(f)L\left(\frac{1}{2},f\right)\omega_{\sym^mf}^z(x)
\\
=
\sum_{q=1}^{+\infty}\frac{1}{\sqrt{q}}e^{-2\pi q/\sqrt{N}}\sum_{n=1}^{+\infty}
\frac{w_m^z(n,q)}{n}e^{-n/x}
+
O(\Err)
\end{multline*}
where
\[%
\Err:= 
N^{-3/8}[\log(2N)]^2\log_2(3N)x^{m/4}[\log(3x)]^{z_m+1}(z_m+m+1)!.
\]
The implicit constant is absolute.
\end{lem}
\begin{proof}
Using \eqref{eq:repd} and Lemma~\ref{lem:fm}, we get
\begin{multline*}
\sum_{f\in\pr}\omegam(f)L\left(\frac{1}{2},f\right)\omega_{\sym^mf}^z(x)
\\
=
\sum_{q=1}^{+\infty}\frac{1}{\sqrt{q}}e^{-2\pi q/\sqrt{N}}\sum_{n=1}^{+\infty}
\frac{w_m^z(n,q)}{n}e^{-n/x}
+O\left(\frac{\tau(N)^2\log_2(3N)}{N^{3/4}}R\right)
\end{multline*}
with
\[%
R
:= 
\sum_{q=1}^{+\infty}\frac{\tau(q)\log(2Nq)}{q^{1/4}}e^{-2\pi q/\sqrt{N}}
\sum_{n=1}^{+\infty}n^{m/4-1}\log(2n)\tau_{(m+1)\abs{z}}(n)e^{-n/x}.
\]
By using 
\[%
\sum_{n\leq t} \frac{\tau_{r}(n)}{n}
\leq [\log(3t)]^{r}
\qquad(t\geq 1, r\geq 1,\,\text{integers}),
\]
we have
\[%
\sum_{n\leq x}\frac{\log(2n)}{n^{1-m/4}}\tau_{(m+1)\abs{z}}(n)e^{-n/x}
\leq
x^{m/4}[\log(3x)]^{z_m+1}
\]
and an integration by parts leads to
\[%
\sum_{n\geq x}\frac{\log(2n)}{n^{1-m/4}}\tau_{(m+1)\abs{z}}(n)e^{-n/x}
\ll_m K
\]
where
\begin{align*}
K
&=
\int_x^{+\infty}\frac{[\log(3t)]^{z_m+1}}{t^{1-m/4}}e^{-t/x}\left(1+\frac{t}{x}\right)\dd t\\
&\leq
x^{m/4}\int_1^{+\infty}[\log(3ux)]^{z_m+1}u^{m/4}e^{-u}(1+1/u)\dd u\\
&\ll_m
x^{m/4}[\log(3x)]^{z_m+1}\int_1^{\infty}u^{m/4+z_m+1}e^{-u}(1+1/u)\dd u\\
&\ll_m
x^{m/4}[\log(3x)]^{z_m}(z_m+m+1)!.
\end{align*}
We conclude with
\[%
\sum_{q=1}^{+\infty}\frac{\tau(q)\log(2Nq)}{q^{1/4}}e^{-2\pi q/\sqrt{N}}
\ll
N^{3/8}[\log(2N)]^2.
\]
\end{proof}

The main term appearing in Lemma~\ref{lem:0918} is studied in the next lemma.
\begin{lem}\label{lem:1819}
Let $m\geq 1$ an integer. There exists $c$ such that, for all $N$ squarefree, $1\leq x^m\leq N^{1/3}$, $z\in \C$,
and $\sigma\in[0,1/3]$ we have
\[%
\sum_{q=1}^{+\infty}\frac{1}{\sqrt{q}}e^{-2\pi q/\sqrt{N}}\sum_{n=1}^{+\infty}
\frac{w_m^z(n,q)}{n}e^{-n/x}
\\=
L^{1,z}\left(\frac{1}{2},1;\st,\sym^m;N\right)
+ O_m(R),
\]
where
\begin{multline*}
R :=  N^{-1/12}e^{c(\abs{z}+1)\log_2(\abs{z}+3)}
\\+ 
x^{-\sigma}\log_2(3N) 
\exp\left\{c(z_m+3)\left[\log_2(z_m+3) 
+ \frac{(z_m+3)^{\sigma/(1-\sigma)}-1}{\sigma\log(z_m+3)}\right]\right\}.
\end{multline*}
The implicit constant depends only on $m$.
\end{lem}
\begin{proof}
Let
\[%
S:= \sum_{q=1}^{+\infty}\frac{1}{\sqrt{q}}e^{-2\pi q/\sqrt{N}}\sum_{n=1}^{+\infty}
\frac{w_m^z(n,q)}{n}e^{-n/x}.
\]
By the definition of $S$, we have $S=S^{>}+S^{\leq}$ with
\begin{multline*}
S^{>}
:= 
\sum_{n_N\mid N^\infty}
\frac{\tau_z(n_N)}{n_N^{m/2+1}}
\sum_{q_N\mid N^\infty}
\frac{\square(n_N^mq_N)}{q_N}
\sum_{\substack{n^{(N)}>x/n_N\\ \left(n^{(N)},N\right)=1}}
\frac{e^{-n_Nn^{(N)}/x}}{n^{(N)}}
\\
\times
\sum_{%
(\nu'_i)_{1\leq i\leq r}\in
\bigtimes{i=1}{r}
[0,m\nu_i]
}
\Biggl\{
\left(
\prod_{j=1}^r
\frac{%
\mu_{\sym^m,\sym^{\nu'_j}}^{z,\nu_j}
}{%
p_j^{\nu'_j/2}
}
\right)
\exp\left(
-\frac{2\pi q_N\prod_{j=1}^{r}p_j^{\nu'_j}}{\sqrt{N}}
\right)
\Biggr\}
\end{multline*}
and
\begin{multline*}
S^{\leq}
:= 
\sum_{n_N\mid N^\infty}
\frac{\tau_z(n_N)}{n_N^{m/2+1}}
\sum_{\substack{n^{(N)}\leq x/n_N\\ \left(n^{(N)},N\right)=1}}
\frac{e^{-n_Nn^{(N)}/x}}{n^{(N)}}
\\
\times%
\sum_{%
(\nu'_i)_{1\leq i\leq r}\in
\bigtimes{i=1}{r}
[0,m\nu_i]
}
\left(
\prod_{j=1}^r
\frac{%
\mu_{\sym^m,\sym^{\nu'_j}}^{z,\nu_j}
}{%
p_j^{\nu'_j/2}
}
\right)
\\
\times
\sum_{q_N\mid N^\infty}
\frac{\square(n_N^mq_N)}{q_N}
\exp\left(
-\frac{2\pi q_N\prod_{j=1}^{r}p_j^{\nu'_j}}{\sqrt{N}}
\right)
\end{multline*}
where $n^{(N)}:= \prod_{j=1}^rp_j^{\nu_j}$.
We have
\begin{equation}\label{eq:msp}
S^>
\ll
R_2
:= 
\sum_{n\mid N^\infty}
\frac{\tau_{\abs{z}}(n)}{n^{m/2+1}}
\sum_{q\mid N^\infty}
\frac{\square(n^mq)}{q}
\sum_{\ell>x/n}
\frac{\widetilde{w}^{z,1/2}_{m,N}(\ell)}{\ell}.
\end{equation}
Moreover, if $n^{(N)}\leq x/n_N$ then
\[%
\prod_{j=1}^rp_j^{\nu'_j}
\leq x^m
\leq N^{1/3}
\]
and
\begin{multline}\label{eq:msm}
\sum_{q_N\mid N^\infty}
\frac{\square(n_N^mq_N)}{q_N}
\exp\left(
-\frac{2\pi q_N\prod_{j=1}^{r}p_j^{\nu'_j}}{\sqrt{N}}
\right)
=%
\\
\sum_{q_N\mid N^\infty}
\frac{\square(n_N^mq_N)}{q_N}
+
O\left(\frac{\tau(n_N^m)}{N^{1/12}}\right).
\end{multline}
Equations \eqref{eq:msp} and \eqref{eq:msm} give $S=P+O(N^{-1/12}R_1+R_2)$ with
\begin{equation*}
P
:= 
\sum_{n_N\mid N^\infty}
\frac{\tau_z(n_N)}{n_N^{m/2+1}}
\sum_{q_N\mid N^\infty}
\frac{\square(n_N^mq_N)}{q_N}
\sum_{\substack{n^{(N)}\leq x/n_N\\ \left(n^{(N)},N\right)=1}}
\frac{\varpi_{m,N}^{z,1/2}(n^{(N)})}{n^{(N)}}
e^{-n^{(N)}/(x/n_N)}
\end{equation*}
and
\[%
R_1
:= 
\sum_{n\mid N^\infty}
\frac
{\tau_{\abs{z}}(n)\tau(n^m)}
{n^{m/2+1}}
\sum_{\ell\leq x/n}\frac{\widetilde{w}^{z,1/2}_{m,N}(\ell)}{\ell}.
\]
Writing
\begin{multline*}
\sum_{\substack{n^{(N)}\leq x/n_N\\ \left(n^{(N)},N\right)=1}}
\frac{\varpi_{m,N}^{z,1/2}(n^{(N)})}{n^{(N)}}
e^{-n^{(N)}/(x/n_N)}
=
W^{z,1/2}_{m,N}(1)
\\
-
\sum_{\substack{\ell> x/n_N\\ \left(\ell,N\right)=1}}
\frac{\varpi_{m,N}^{z,1/2}(\ell)}{\ell}
+
\sum_{\substack{\ell\leq x/n_N\\ \left(\ell,N\right)=1}}
\frac{\varpi_{m,N}^{z,1/2}(\ell)}{\ell}
\left[e^{-\ell/(x/n_N)}-1\right]
\end{multline*}
we get, by Lemma~\ref{lem:j46},
\[%
P= L^{1,z}\left(\frac{1}{2},1;\st,\sym^m;N\right)+ O(R_2 + R_3)
\]
with
\begin{equation*}
R_3:= 
\sum_{n\mid N^\infty}\frac{\tau_{\abs{z}}(n)}{n^{m/2+1}}
\sum_{q_\mid N^\infty}
\frac{\square(n^mq)}{q}
\sum_{\substack{\ell\leq x/n\\ \left(\ell,N\right)=1}}
\frac{\widetilde{w}_{m,N}^{z,1/2}(\ell)}{\ell}
\left[1-e^{-\ell/(x/n)}\right].
\end{equation*}
Lemma~\ref{lem:j45} gives
\[%
R_1\ll\exp\left[
c(z_m+3)\log_2(z_m+3)
\right].
\]
We have
\begin{align*}
R_3
&\ll
\sum_{n\mid N^\infty}\frac{\tau_{\abs{z}}(n)}{n^{m/2+1}}
\sum_{q\mid N^\infty}
\frac{\square(n^mq)}{q}
\sum_{\ell\leq x/n}
\frac{\widetilde{w}_{m,N}^{z,1/2}(\ell)}{\ell}
\cdot
\frac{\ell n}{x}\\
&\ll
x^{-\sigma}
\sum_{n\mid N^\infty}\frac{\tau_{\abs{z}}(n)}{n^{m/2+1-\sigma}}
\sum_{q_\mid N^\infty}
\frac{\square(n^mq)}{q}
\sum_{\ell=1}^{+\infty}
\frac{\widetilde{w}_{m,N}^{z,1/2}(\ell)}{\ell^{1-\sigma}}
\end{align*}
for all $\sigma\in[0,1/2[$ and Lemma~\ref{lem:j45} gives
\begin{multline*}%
R_3
\ll%
\\
x^{-\sigma}\log_2(3N)
\exp\left\{c(z_m+3)\left[\log_2(z_m+3) 
+ \frac{(z_m+3)^{\sigma/(1-\sigma)}-1}{\sigma\log(z_m+3)}\right]\right\}.
\end{multline*}
Next, for all $\sigma\in[0,1/2[$ , Rankin's method and Lemma~\ref{lem:j45} give
\begin{multline*}%
R_2
\ll%
\\
x^{-\sigma}\log_2(3N)
\exp\left\{c(z_m+3)\left[\log_2(z_m+3) 
+ \frac{(z_m+3)^{\sigma/(1-\sigma)}-1}{\sigma\log(z_m+3)}\right]\right\}.
\end{multline*}
\end{proof} 

Next, given $\eta\in]0,1/100[$, denote by $\prp$\label{page:td} the subset of $\pr$ consisting of forms $f$ such that
$L(s,\sym^mf)$ has no zeros in the half strip
\[%
\Re s\geq 1-4\eta \qquad \abs{\Im s}\leq 2[\log(2N)]^3 
\]
and $\prm$ the complementary subset. By \cite[Proposition 5.3]{CoMi04}, for all $m\geq 1$, there exists
$\xi>0$ and $A>0$ (both depending on $m$) such that for all $\eta\in]0,1/100[$ and squarefree $N$ we have
\begin{equation*}
\#\prm\leq \xi N^{A\eta}[\log(2N)]^{\xi}.
\end{equation*}
By \cite[Lemmas 4.1 and 4.2]{CoMi04} 
there exists, for all $m\geq 1$, a constant $B$ (depending on $m$) such that, for all $z\in\C$ and $f\in\prm$,
we have
\begin{equation}\label{eq:majlsy}
L(1,\sym^mf)^z\ll_m [\log(2N)]^{B\abs{\Re z}}
\end{equation}
Using the convexity bound (see \cite[Lecture 4]{mic03} for better bounds that we do not need here)
\[%
L\left(\frac{1}{2},f\right)\ll N^{1/4}
\]
and
\[%
\omegam(f)=\frac{\pi^2}{\varphi(N)L(1,\sym^2f)}\ll\frac{\log(2N)\log_2(3N)}{N}
\]
and by  \eqref{eq:majlsy} we get
\[%
\sum_{f\in\prm}\omegam(f)L\left(\frac{1}{2},f\right)L(1,\sym^mf)^z
\ll_m
N^{A\eta-3/4}[\log(2N)]^{B\abs{\Re z}+C},
\]
$A$, $B$ and $C$ being constants depending only on $m$ 
so that
\begin{multline*}
\sum_{f\in\pr}\omegam(f)L\left(\frac{1}{2},f\right)L(1,\sym^mf)^z
\\=
\sum_{f\in\prp}\omegam(f)L\left(\frac{1}{2},f\right)L(1,\sym^mf)^z
\\+O_m\left(
N^{A\eta-3/4}[\log(2N)]^{B\abs{\Re z}+C}
\right).
\end{multline*}
Next, there exists a constant $D>0$, depending only on $m$, such that
\[%
L(1, \sym^{m}f)^z
= \omega_{\sym^{m}f}^z(x) + O(R_1),
\]
with
\[%
R_1 :=  x^{-1/\log_2(3N)} e^{D\abs{z}\log_3(20N)} [\log(2N)]^3
+ e^{D\abs{z}\log_2(3N)-[\log(2N)]^2}
\]
(see \cite[Proposition 5.6]{CoMi04}) and, since by positivity (see \cite{Guo96} and \cite{FH95}) 
and Lemma~\ref{lem:moyLdem} we have
\[%
\sum_{f\in\prp}\omegam(f)L\left(\frac{1}{2},f\right)\ll 1,
\] 
we obtain
\begin{multline*}
\sum_{f\in\pr}\omegam(f)L\left(\frac{1}{2},f\right)L(1,\sym^mf)^z
\\=
\sum_{f\in\prp}\omegam(f)L\left(\frac{1}{2},f\right)\omega_{\sym^{m}f}^z(x)+O_m(R_2)
\end{multline*}
with
\[%
R_2:= R_1+N^{A\eta-3/4}[\log(2N)]^{B\abs{\Re z}+C}.
\]
Now, since $\abs{\omega_{\sym^m f}^z(x)}\leq \iota(\epsilon)^{\abs{\Re z}}x^\epsilon$, 
where $\iota(\epsilon)>1$ depends on $\epsilon$ and $m$, we reintroduce 
the forms of $\prm$ obtaining
\begin{multline*}
\sum_{f\in\pr}\omegam(f)L\left(\frac{1}{2},f\right)L(1,\sym^mf)^z
\\=
\sum_{f\in\pr}\omegam(f)L\left(\frac{1}{2},f\right)\omega_{\sym^{m}f}^z(x)+O_m(R_3)
\end{multline*}
with
\begin{multline*}
R_3:= 
x^{-1/\log_2(3N)} e^{D\abs{z}\log_3(20N)}[\log(2N)]^3
\\
+x^\epsilon N^{A\eta-3/4}[\iota(\epsilon)\log(2N)]^{B\abs{\Re z}+C}
+e^{D\abs{z}\log_2(3N)-[\log(2N)]^2}.
\end{multline*}
Lemmas~\ref{lem:0918} and~\ref{lem:1819} with
$\eta=\epsilon=1/(100m)$, $x^m=N^{1/10}$ and \[\sigma=c'(m)/\log(\abs{z}+3)\] with $c'(m)$ large enough and depending
on $m$ leads to Theorem~\ref{thmint:A}.

\subsection{\proofname{} of  Proposition~\ref{propint:B}}

For the proof of Proposition~\ref{propint:B}, we write
\begin{multline*}
L^{1,z}\left(\frac{1}{2},1;\st,\sym^m;N\right)
=
L^{1,z}\left(\frac{1}{2},1;\st,\sym^m\right)
\\
\times
\LB^z_m(N)\prod_{p\mid N}
\left(
\int_{\su}D(p^{-1/2},\st,g)D(p^{-1},\sym^m,g)^z\dd g
\right)^{-1}.
\end{multline*}
We use
\[%
\LB^z_m(N)=%
1+O\left(\frac{(\abs{z}+1)\omega(N)}{P^-(N)^{\min\{m/2+1,2\}}}\right)%
\]
which is uniform for all $z$ and $N$ such that 
\[(\abs{z}+1)\omega(N)\leq P^-(N)^{\min\{m/2+1,2\}}\]
and Lemma~\ref{lem:edf}
to get
\[%
L^{1,z}\left(\frac{1}{2},1;\st,\sym^m;N\right)%
=%
L^{1,z}\left(\frac{1}{2},1;\st,\sym^m\right)%
[1+O_m(\Err)]%
\]
where
\[%
\Err:= 
\frac{\omega(N)}{P^-(N)}
+
\frac{(\abs{z}+1)\omega(N)}{P^-(N)^{3/2}}
+
\frac{(\abs{z}+1)^2\omega(N)}{P^-(N)^{2}}
\]
uniformely for
\begin{equation*}
\begin{cases}
N\in\Ncri{%
\max\bigl\{
\omega(\cdot)^{1/2},\,
[(\abs{z}+1)\omega(\cdot)]^{2/3},\, 
[\abs{z}^2\omega(\cdot)]^{1/2}
\bigl\}
}
,\text{}\\
z\in\C.
\end{cases}
\end{equation*}
 \section{Twisting by $L(1,\sym^2f)$}
In this section, we sketch the proofs of Theorem~\ref{thmint:B} and Proposition~\ref{propint:D}.
The proof of Theorem~\ref{thmint:B} is very similar to the one of Theorem~\ref{thmint:A}.

Let $z\in\C$ and $x\geq 1$, define
\begin{equation}\label{eq:auomeg}
\omega_{\sym^2f,\sym^mf}^{1,z}(x):= \sum_{n=1}^{+\infty}\frac{\lambda_{\sym^2f,\sym^mf}^{1,z}(n)}{n}e^{-n/x}.
\end{equation}
for all $f\in\pr$ and obtains the
\begin{lem}
Let $N$ be a squarefree integer, $m\in\Z_{>0}$, $x\geq 1$ and $z\in\C$. 
Then
\[%
\sum_{f\in\pr}\omegam(f)\omega_{\sym^2f,\sym^mf}^{1,z}(x)
=
\frac{\varphi(N)}{N}
\sum_{n=1}^{+\infty}
\frac{w_{2,m}^{1,z}(n)}{n}e^{-n/x}
+O(\Err)
\]
with
\[%
\Err:= 
\frac{\tau(N)^2\log(2N)\log_2(3N)}{N}x^{m/4}(\log 3x)^{z_m+3}(z_m+m+4)!.
\]
The implicit constant is absolute and $w_{2,m}^{1,z}(n)$ has been defined in Lemma~\ref{lem:siegfried}.
\end{lem}

Next, we have the
\begin{lem}
Let $m\geq 1$ an integer. There exists $c$ such that, for all $N$ squarefree, $1\leq x^m\leq N^{1/3}$, $z\in \C$,
and $\sigma\in[0,1/{3m}]$ we have
\[%
\sum_{n=1}^{+\infty}
\frac{w_{2,m}^{1,z}(n)}{n}e^{-n/x}=
L^{1,z}\left(1,1;\sym^2,\sym^m;N\right)
+ O_m(R),
\]
where
\[%
R :=  
\\
\frac{\log_2(3N)}{x^\sigma} 
\exp\left\{c(z_m+3)\left(\log_2(z_m+3) 
+ \frac{(z_m+3)^{\sigma/(1-\sigma)}-1}{\sigma\log(z_m+3)}\right)\right\}.
\]
The implicit constant depends only on $m$.
\end{lem} 

The conclusion of the proof of Theorem~\ref{thmint:B} is the same as the one of Theorem~\ref{thmint:A} after 
having introduced the exceptional set
\[%
\prd:= \pr\setminus\left(\prpdd{2}\cap\prp\right).
\] 
The proof of Proposition~\ref{propint:D} follows from Lemma~\ref{lem:edf} in the same way as Proposition~\ref{propint:B}.
\section{Asymptotic of the moments}
\subsection{\proofname{} of Proposition~\ref{propint:asymp}}

We give the proof for
$L^{1,\pm r}\left(\frac{1}{2},1;\st,\sym^m\right)$ since the method is similar in the two cases.

Write
\[%
\psi_{m,1}^{\pm r}(p):= \int_{\su}D(p^{-1/2},\st,g)D(p^{-1},\sym^m,g)^{\pm r}\dd g.
\]
By Lemma~\ref{lem:diable}, we have
\[%
\sum_{p\geq (m+1)r+3}\log\psi_{m,1}^{\pm r}(p)\ll_m\frac{r}{\log r}.
\]
By \eqref{eq:lisyn} we get
\[%
\left(1+\frac{1}{\sqrt{p}}\right)^{-2}D(p^{-1},\sym^m,g)
\leq
\psi_{m,1}^{\pm r}(p)
\leq
\left(1-\frac{1}{\sqrt{p}}\right)^{-2}D(p^{-1},\sym^m,g)
\]
and then
\begin{equation}\label{eq:1642}
\sum_{p\leq (m+1)r+3}\log\psi_{m,1}^{\pm r}(p)
=
\sum_{p\leq (m+1)r+3}\log\Upsilon_{m,1}^{\pm r}(p)
+
O_m\left(\sqrt{r}\log_2(3r)\right)
\end{equation}
with
\[%
\Upsilon_{m,1}^{\pm r}(p):= \int_{\su}D(p^{-1},\sym^m,g)^{\pm r}\dd g.
\]
The right hand side of \eqref{eq:1642} has been evaluated in \cite[\S 2.2.1]{CoMi04} and was founded to be
\[%
\sym^m_{\pm}r\log_2 r
+\sym^{m,1}_{\pm}r
+O_m\left(\frac{r}{\log r}\right)
\]
which ends the proof.

\subsection{\proofname{} of Corollary~\ref{lem:dblcdn}}

Let $r\geq 0$.
Define
\[%
\Theta(N) 
:=  
\sum_{g\in\pr}
\omega(g)
L\left(\frac{1}{2},g\right)
\qquad \text{and}\qquad
\Omega(f) 
:=  
\frac{
\omega(f) L\left(\frac{1}{2}, f\right)
}
{\Theta(N)}.
\]
For $N\in\Ncri{\log^{1/2}}$, we have
\[%
\Theta(N)\sim 1 \qquad (N\to +\infty)
\]
(see Lemma~\ref{lem:moyLdem}).
Since $L\left(\frac{1}{2},f\right)\geq 0$, by Theorem~\ref{thmint:A}, and Propositions~\ref{propint:B} 
and~\ref{propint:asymp} we get 
\begin{align*}
\sum_{\substack{%
f\in\pr\\
L\left(\frac{1}{2},f\right)>0
}}
\Omega(f) 
L(1,\sym^mf)^{r}
& = 
\frac{1}{\Theta(N)}
\sum_{f\in\pr} 
\omega(f) 
L\left(\frac{1}{2}, f\right) L(1,\sym^mf)^{r}
\\
& = 
[1+o(1)]
e^{%
\sym^m_+r
\log\left\{
[1+o(1)]\exp\left(\frac{\sym^{m,1}_+}{\sym^m_+}\right)\log r
\right\}
}
\end{align*}
uniformly for all $r\leq c\log N/\log_2(3N)\log_3(20N)$.
Since
\[%
\sum_{\substack{%
f\in\pr\\
L\left(\frac{1}{2},f\right)>0
}}
\Omega(f) 
=
\sum_{f\in\pr}\Omega(f)=1
\]
we obtain, by positivity,a function
$f\in\pr$ such that
\[%
L(1,{\rm sym}^mf)^{r}
\geq\{1+o(1)\} 
e^{
\sym^m_+r\log\left\{
[1+o(1)]\exp(\sym^{m,1}_+/\sym^m_+)\log r
\right\}
}
\]
and $L\left(\frac{1}{2},f\right)>0$.

We obtain the announced minoration with 
$r= c\log N/(\log_2(3N))^2$.
The majoration is obtained in the same way, taking the negative moments.
\section{Hecke eigenvalues}
\subsection{\proofname{} of Proposition~\ref{propint:rubin}}

Following step by step the proof given by Granville \& Soundararajan in the case of Dirichlet
characters \cite[Lemma 8.2]{gra01c}, we get under Grand Riemann Hypothesis
\[%
\log
L(1,\sym^mf)=
\sum_{2\leq n\leq\log^2(2N)\log_2^4(3N)}\frac{\Lambda_{\sym^mf}(n)}{n\log n}
+
O_m\left(1\right)
\]
where $\Lambda_{\sym^m}(n)$ is the function defined by
\[%
-\frac{L'(s,\sym^mf)}{L(s,\sym^mf)}=: \sum_{n=1}^{+\infty}\frac{\Lambda_{\sym^m}(n)}{n^s}\qquad (\Re s>1)
\]
that is
\[%
\Lambda_{\sym^m}(n)
=
\begin{cases}
\chi_{\sym^m}[g(\theta_{f,p})^\nu]\log p & \text{if $n=p^\nu$ with $p\nmid N$}\\
\lambda_f(p)^{m\nu}\log p & \text{if $n=p^\nu$ with $p\mid N$}\\ 
0 & \text{otherwise.}
\end{cases}
\]
If $\nu>1$, then
\[%
\abs{\frac{\Lambda_{\sym^mf}(p^{\nu})}{p^{\nu}\log(p^{\nu})}}
\leq
\frac{m+1}{p^{\nu}}
\]
hence
\[%
\log L(1,\sym^mf)=\sum_{p\leq\log^2(2N)\log_2^4(3N)}
\frac{\Lambda_{\sym^mf}(p)}{p\log p}
+O(1).
\]
From $\Lambda_{\sym^mf}(p)=\lambda_f(p^m)\log p$ we deduce
\begin{equation*}
\log L(1,\sym^mf)=\sum_{p\leq\log^2(2N)\log_2^4(3N)}
\frac{\lambda_{f}(p^m)}{p}
+O(1).
\end{equation*}
Since
\begin{align*}
\sum_{\log(2N)\leq p\leq \log^2(2N)\log_2^4(3N)}\frac{\abs{\lambda_{f}(p^m)}}{p}
& \leq
(m+1)\sum_{\log(2N)\leq p\leq \log^2(2N)\log_2^4(3N)}\frac{1}{p}\\
& \ll_m1
\end{align*}
we get
\begin{equation}\label{eq:jacob}
\log L(1,\sym^mf)=\sum_{p\leq\log(2N)}
\frac{\lambda_{f}(p^m)}{p}
+O_m(1).
\end{equation}
Let $N\in\Ncri{\log^{3/2}}$ and $f\in\prgde{m}$, equation \eqref{eq:jacob} then leads to
\[%
\sum_{p\leq\log(2N)}\frac{\lambda_f(p^m)}{p}\geq\sym^m_+\log_3(20N)+O_m(1)
\]
and we deduce
\[%
\sum_{p\leq\log(2N)}\frac{\sym^m_+-\lambda_f(p^m)}{p}\ll_m1.
\]
For $\xi(N)\leq\log_3(20N)$, we get
\begin{align*}
\sum_{\substack{%
p\leq\log(2N)
\\
\lambda_{f}(p^m)\geq\sym^m_+-\xi(N)/\log_3(20N)
}} 
\frac{1}{p}
&=
\sum_{%
p\leq\log(2N)
} 
\frac{1}{p}\\
& \phantom{= \log_3(20N)}
-
\sum_{\substack{%
p\leq\log(2N)
\\
\lambda_{f}(p^m)<\sym^m_+-\xi(N)/\log_3(20N)
}} 
\frac{1}{p}
\\
& = \log_3(20N) 
\left\{
1 + O_{\epsilon, m}\left(\frac{1}{\xi(N)}\right)
\right\}.
\end{align*}
We conclude by using
\[%
\sum_{\log^{\epsilon}(3N)<p<\log(2N)}\frac{1}{p}\ll 1.
\]

\subsection{\proofname{} of Proposition~\ref{propint:tram}}

Let $N\in\Ncri{\log^{3/2}}$.
Taking $m=2$ in \eqref{eq:jacob} gives
\[%
\sum_{p\leq\log(2N)}\frac{\lambda_f(p^2)}{p}+O(1)
=
\log L(1,\sym^2f).
\]
Since $\sym^2_-=1$, if $f\in\prpte{2}$, we deduce
\[%
\sum_{p\leq\log(2N)}\frac{\lambda_f(p^2)}{p}
\leq
-\log_3(20N)+O(1).
\]
If $p\mid N$, then $\lambda_f(p^2)=\lambda_f(p)^2$ and
\[%
\sum_{\substack{p\leq\log(2N)\\ p\mid N}}\frac{1}{p}
=O(1);
\]
if $p\nmid N$, then $\lambda_f(p^2)=\lambda_f(p)^2-1$. We thus have
\[%
\sum_{p\leq\log(2N)}\frac{\lambda_f(p)^2-1}{p}
\leq
-\log_3(20N)+O(1)
\]
hence
\begin{equation}\label{eq:cafe}
\sum_{p\leq\log(2N)}\frac{\lambda_f(p)^2}{p}
\ll 1.
\end{equation}
For $\xi(N)\leq\log_3(20N)$, we deduce
\[%
\sum_{\substack{p\leq\log(2N)\\\abs{\lambda_f(p)}\geq[\xi(N)/\log_3(20N)]^{1/2}}}\frac{\lambda_f(p)^2}{p}
\ll
\frac{\log_3(20N)}{\xi(N)}
\]
which leads to the announced result.

\section{Simultaneous extremal values}

\subsection{\proofname{} of Proposition~\ref{propint:simul}}

Prove the first point.
Let $C>0$, $N\in\Ncri{\log}$ and $f\in\pr$ such that
\[%
L(1,\sym^{2}f)
\leq C
\left[\log_2(3N)\right]^{-\sym^2_{-}}
\]
and
\[%
L(1,\sym^{4}f)
\leq C
\left[\log_2(3N)\right]^{-\sym^4_{-}}.
\]
Equation \eqref{eq:jacob} with $m=4$ gives
\[%
\sum_{\substack{p\leq\log(2N)\\ p\nmid N}}\frac{\lambda_f(p^4)}{p}+O(1)
\leq 
-\sym^4_-\log_3(20N)
\]
since the contribution of $p$ dividing $N$ is bounded (using \eqref{eq:krivine}). Expanding $\lambda_f(p^4)$ thanks to
\eqref{eq:multhecke} we deduce
\[%
\sum_{\substack{p\leq\log(2N)\\ p\nmid N}}\frac{\lambda_f(p)^4-3\lambda_f(p)^2+1}{p}+O(1)
\leq 
-\sym^4_-\log_3(20N).
\]
Reinserting \eqref{eq:cafe} (again, we remove easily the contribution of $p$ dividing $N$), we are led to
\[%
\sum_{\substack{p\leq\log(2N)\\ p\nmid N}}\frac{\lambda_f(p)^4+1}{p}
\leq 
-\sym^4_-\log_3(20N)+O(1).
\]
The right hand side tends to $-\infty$ while the left one is positive, so we get a contradiction.

Prove next the second point. Assume that
\[%
L(1,\sym^{2}f)
\geq C
\left[\log_2(3N)\right]^{\sym^2_{+}}.
\]
By Cauchy-Schwarz inequality and \eqref{eq:jacob}, we have
\begin{equation}\label{eq:fluid}
(\sym^2_+)^2[\log_2(3N)+O(1)]\leq
\sum_{\substack{p\leq\log(2N)\\ p\nmid N}}\frac{\lambda_f(p^2)^2}{p}.
\end{equation}
Further, from $X_4=X_2^2-X_2-1$, we deduce
\[%
\sum_{\substack{p\leq\log(2N)\\ p\nmid N}}\frac{\lambda_f(p^4)}{p}
=
\sum_{\substack{p\leq\log(2N)\\ p\nmid N}}\frac{\lambda_f(p^2)^2-\lambda_f(p^2)-1}{p}
\]
and \eqref{eq:fluid} and $\abs{\lambda_f(p^2)}\leq\sym^2_+$ imply
\[%
\sum_{\substack{p\leq\log(2N)\\ p\nmid N}}\frac{\lambda_f(p^4)}{p}
\geq
[(\sym^2_+)^2-\sym_+^2-1]\log_3(20N)+O(1)
\]
which leads to the result by \eqref{eq:jacob} since
\[%
(\sym^2_+)^2-\sym_+^2-1=\sym^4_+.
\]

\subsection{\proofname{} of Proposition~\ref{propint:simulde}}

From 
\[%
X_m^2=\sum_{j=2}^mX_{2j}+X^2
\]
we deduce
\begin{align*}
\sum_{\substack{p\leq\log(2N)\\ p\nmid N}}\frac{\lambda_f(p^m)^2}{p}
&=
\sum_{\substack{p\leq\log(2N)\\ p\nmid N}}\sum_{j=2}^m\frac{\lambda_f(p^{2j})}{p}
+
\sum_{\substack{p\leq\log(2N)\\ p\nmid N}}\frac{\lambda_f(p)^2}{p}\\
&\leq
(m+3)(m+1)\log_3(20N)+O(1)
\end{align*}
by \eqref{eq:cafe} and $\abs{\lambda_f(p^{2j})}\leq 2j+1$. Furthermore
\begin{align*}
[\sym^m_+\log_3(20N)]^2
&=
\left(
\sum_{\substack{p\leq\log(2N)\\ p\nmid N}}\frac{\lambda_f(p^m)}{p}
\right)^2\\
&\leq
[\log_3(20N)+O(1)]
\sum_{\substack{p\leq\log(2N)\\ p\nmid N}}\frac{\lambda_f(p^m)^2}{p}
\end{align*}
so that
\[%
(\sym^m_+)^2\leq(m+3)(m-1)
\]
which contradicts $\sym^m_+=(m+1)^2$.
\section{An index of notation}\label{sec:index}

$%
\begin{array}{|l|l||l|l||l|l|}
\hline%
\gamma^* \! & \eqref{eq:gamet} \! & %
\lambda_{\sym^mf}^{z,\nu}(\phantom{n}) \! & \eqref{eq:devD} \! & %
\chi \! & \! \text{p. \pageref{def:char}}\!\! \\
\delta(\phantom{m},\phantom{n}) \! & \text{\S~\ref{pa:us}} \! & %
\lambda_{\sym^2,\sym^m}^{1,z}(\phantom{n}) \! & \eqref{eq:conf} \! & %
\omegam \! & \! \eqref{eq:hawe}\!\! \\
\Delta(\phantom{a},\phantom{b},\phantom{c}) \! & \eqref{eq:valDel} \! & %
\lambda_{\sym^2,\sym^m}^{1,z,\nu}(\phantom{n}) \! & \eqref{eq:dvu} \! & %
\omega^{z}_{\sym^mf}(x) \! & \! \eqref{eq:omega}\!\! \\
\Delta(\phantom{a},\phantom{b},\phantom{c};\phantom{d},\phantom{e}) \! & \eqref{eq:monpey} \! & %
\mu_{\sym^m,\sym^{m'}}^{z,\nu} \! &  \eqref{eq:ei} \! & %
\omega^{1,z}_{\sym^2f,\sym^mf}(\phantom{n}) \! & \! \eqref{eq:auomeg}\!\! \\
\epsilon_f(N) \! & \eqref{eq:tirage} \! & %
\mu_{\sym^2,\sym^m,\sym^{m'}}^{1,z,\nu} \! & \eqref{eq:nsle} \! & %
\varpi_{m,N}^{z,\rho}(\phantom{n}) \! & \! \eqref{eq:chri2} \!\! \\
\zeta^{(N)} \! & \eqref{eq:caN} %
\! & \rho \! & \text{\S~\ref{pa:us}} \! & %
\widetilde{w}_{m,N}^{z,\rho}(\phantom{n}) \! & \! \eqref{eq:chri3} \!\! \\
\lambda_f(\phantom{n}) \! & \eqref{eq:deff} \! & %
\sigma \! & \text{\S~\ref{pa:us}} \! & %
\! & \!\! \\
\lambda_{\sym^mf}^{z}(\phantom{n}) \! & \eqref{eq:seco} \! & %
\tau_z(\phantom{n}) \! & \eqref{eq:zz}\! & %
\! & \!\! \\ %
\hline%
\end{array}
$

$%
\begin{array}{|l|l||l|l|}
\hline%
D(\phantom{a},\phantom{b},\phantom{c}) & \eqref{eq:fato} &%
n^{(N)} & \text{\S~\ref{pa:us}}\\
g(\phantom{a}) & \eqref{eq:defgt} &%
\Ncri{\phantom{g}} & \eqref{eq:defng}\\
\pr & \text{p. \pageref{page:un}} &%
P^-(\phantom{N}) & \text{p. \pageref{page:cin}}\\
\prp & \text{p. \pageref{page:td}} &%
\sym^m_{\pm} & \eqref{eq:defsy}\\
\prm & \text{p. \pageref{page:td}} &%
\sym^{m,1}_{\pm} & \eqref{eq:defsyu}\\
\prgde{m} & \eqref{eq:dhp} &%
w_m^z(\phantom{a},\phantom{b}) & \eqref{eq:defterpr} \\
\ell(m,\vec{\nu}) & \eqref{eq:ellem} &%
w^{1,z}_{2,m}(\phantom{a}) & \text{Lemme~\ref{lem:siegfried}} \\
\ell(2,m;\vec{\nu},\vec{\nu'}) & \eqref{eq:miap} &%
W^{z,\rho}_{m,N}(\phantom{a}) & \eqref{eq:chri} \\
L^{1,z}\left(\frac{1}{2},1;\st,\sym^m;N\right) & \eqref{eq:eli1} &%
X_m & \eqref{eq:trtch} \\
L^{1,z}\left(\frac{1}{2},1;\st,\sym^m\right) & \eqref{eq:eli2} &%
\LB_m^z(N) & \eqref{eq:dtm} \\
L^{1,z}\left(1,1;\sym^2,\sym^m,N\right) & \eqref{eq:heure} &%
\LB_{2,m}^{1,z}(N) & \eqref{eq:lenti} \\
L^{1,z}\left(1,1;\sym^2,\sym^m\right) & \eqref{eq:charo} &%
z_m & \eqref{eq:defzm} \\
n_N & \text{\S~\ref{pa:us}} &%
    &      \\
\hline%
\end{array}
$%

$%
\begin{array}{|l|l|}
\hline%
\square & \text{\S~\ref{pa:us}}\\
\square_N(\phantom{a}) & \eqref{eq:caN}\\
\un_N & \text{\S~\ref{pa:us}}\\
\un^{(N)} & \text{\S~\ref{pa:us}}\\
\hline%
\end{array}
$



\begin{thebibliography}{xx}
\bibitem[CFKRS03]{CFKRS03}
J.~B. Conrey, D.~W. Farmer, J.~P. Keating, M.~O. Rubinstein, and N.~C. Snaith,
  \emph{Autocorrelation of random matrix polynomials}, Comm. Math. Phys.
  \textbf{237} (2003), no.~3, 365--395.

\bibitem[CFKRS05]{CFKRS05}
J.~B. Conrey, D.~W. Farmer, J.~P. Keating, M.~O. Rubinstein, and N.~C. Snaith,
 \emph{Integral moments of {$L$}-functions}, Proc. London Math. Soc.
  (3) \textbf{91} (2005), no.~1, 33--104. 

\bibitem[CM04]{CoMi04}
J.~Cogdell and P.~Michel, \emph{On the complex moments of symmetric power
  {$L$}-functions at {$s=1$}}, Int. Math. Res. Not. (2004), no.~31, 1561--1617.
  

\bibitem[Eic54]{eic54}
M.~Eichler, \emph{{Quatern{\"a}re quadratische Formen und die Riemannsche
  Vermutung}}, Archiv der Mathematik \textbf{V} (1954), 355--366.

\bibitem[Ell73]{ell73}
P.~D. T.~A. Elliott, \emph{On the distribution of the values of quadratic
  {$L$}-series in the half-plane {$\sigma >\frac{1}{2}$}}, Invent. Math.
  \textbf{21} (1973), 319--338. 

\bibitem[FH95]{FH95}
S.~Friedberg and J.~Hoffstein, \emph{Nonvanishing theorems for automorphic
  {$L$}-functions on {${\rm GL}(2)$}}, Ann. of Math. (2) \textbf{142} (1995),
  no.~2, 385--423. 

\bibitem[FOP04]{MR2114776}
Sharon Frechette, Ken Ono, and Matthew Papanikolas, \emph{Combinatorics of
  traces of {H}ecke operators}, Proc. Natl. Acad. Sci. USA \textbf{101} (2004),
  no.~49, 17016--17020 (electronic). 

\bibitem[GHL94]{ghl94}
D.~Goldfeld, J.~Hoffstein, and D.~Lieman, \emph{An effective zero-free region},
  Ann. of Math. (2) \textbf{140} (1994), no.~1, 177--181, Appendix of
  \cite{hol94}.

\bibitem[GJ78]{gej78}
Stephen Gelbart and Herv{\'e} Jacquet, \emph{A relation between automorphic
  representations of {${\rm GL}(2)$} and {${\rm GL}(3)$}}, %
Ann. Sci. \'Ecole Norm. Sup. (4) \textbf{11} (1978), no.~4, 471--542. 

\bibitem[GS01]{gra01c}
A.~Granville and K.~Soundararajan, \emph{Large character sums}, J. Amer. Math.
  Soc. \textbf{14} (2001), no.~2, 365--397 (electronic). 

\bibitem[GS03]{GS03}
A.~Granville and K.~Soundararajan,
\emph{The distribution of values of {$L(1,\chi\sb d)$}}, Geom. Funct.
  Anal. \textbf{13} (2003), no.~5, 992--1028. 

\bibitem[Guo96]{Guo96}
J.~Guo, \emph{On the positivity of the central critical values of automorphic
  {$L$}-functions for {${\rm GL}(2)$}}, Duke Math. J. \textbf{83} (1996),
  no.~1, 157--190. 

\bibitem[HL94]{hol94}
J.~Hoffstein and P.~Lockhart, \emph{Coefficients of {M}aass forms and the
  {S}iegel zero}, Ann. of Math. (2) \textbf{140} (1994), no.~1, 161--181, With
  an appendix by D. Goldfeld, J. Hoffstein and D. Lieman. 

\bibitem[HR95]{hor95}
J.~Hoffstein and D.~Ramakrishnan, \emph{Siegel zeros and cusp forms}, Internat.
  Math. Res. Notices (1995), no.~6, 279--308. 

\bibitem[HR04]{MR2139690}
Laurent Habsieger and Emmanuel Royer, \emph{{$L$}-functions of automorphic
  forms and combinatorics: {D}yck paths}, Ann. Inst. Fourier (Grenoble)
  \textbf{54} (2004), no.~7, 2105--2141 (2005). 

\bibitem[Igu59]{igu59}
Jun-ichi Igusa, \emph{Kroneckerian model of fields of elliptic modular
  functions}, Amer. J. Math. \textbf{81} (1959), 561--577. 

\bibitem[IK04]{IwKo04}
Henryk Iwaniec and Emmanuel Kowalski, \emph{Analytic number theory}, American
  Mathematical Society Colloquium Publications, vol.~53, American Mathematical
  Society, Providence, RI, 2004. 

\bibitem[ILS00]{ils00}
H.~Iwaniec, W.~Luo, and P.~Sarnak, \emph{Low lying zeros of families of
  ${L}$-functions}, Inst. Hautes \'Etudes Sci. Publ. Math. (2000), no.~91,
  55--131 (2001).

\bibitem[IS00]{IS00}
H.~Iwaniec and P.~Sarnak, \emph{The non-vanishing of central values of
  automorphic {$L$}-functions and {L}andau-{S}iegel zeros}, Israel J. Math.
  \textbf{120} (2000), part A, 155--177. 

\bibitem[Kim03]{kim03}
Henry~H. Kim, \emph{Functoriality for the exterior square of {${\rm GL}\sb 4$}
  and the symmetric fourth of {${\rm GL}\sb 2$}}, J. Amer. Math. Soc.
  \textbf{16} (2003), no.~1, 139--183 (electronic), With appendix 1 by Dinakar
  Ramakrishnan and appendix 2 by Kim and Peter Sarnak. 
  

\bibitem[KMV00]{KMV00}
E.~Kowalski, P.~Michel, and J.~VanderKam, \emph{Mollification of the fourth
  moment of automorphic {$L$}-functions and arithmetic applications}, Invent.
  Math. \textbf{142} (2000), no.~1, 95--151. 

\bibitem[KS02a]{KS02b}
Henry~H. Kim and Freydoon Shahidi, \emph{Functorial products for {${\rm GL}\sb
  2\times{\rm GL}\sb 3$} and the symmetric cube for {${\rm GL}\sb 2$}}, Ann. of
  Math. (2) \textbf{155} (2002), no.~3, 837--893, With an appendix by Colin J.
  Bushnell and Guy Henniart. 

\bibitem[KS02b]{KS02}
H.H. Kim and F.~Shahidi, \emph{Cuspidality of symmetric powers with
  applications}, Duke Math. J. \textbf{112} (2002), no.~1, 177--197. 

\bibitem[Lit28]{L28}
J.E. Littlewood, \emph{{On the class number of the corpus $P(\sqrt{-k})$}},
  Proc. London Math. Soc. \textbf{27} (1928), 358--372.

\bibitem[Luo99]{Luo99}
W.~Luo, \emph{Values of symmetric square {$L$}-functions at {$1$}}, J. Reine
  Angew. Math. \textbf{506} (1999), 215--235. 

\bibitem[Luo01]{Luo00}
W.~Luo, \emph{Nonvanishing of {$L$}-values and the {W}eyl law}, Ann. of Math.
  (2) \textbf{154} (2001), no.~2, 477--502. 

\bibitem[LW06]{LW04}
Yuk-Kam Lau and Jie Wu, \emph{A density theorem on automorphic {$L$}-functions
  and some applications}, Trans. Amer. Math. Soc. \textbf{358} (2006), no.~1,
  441--472 (electronic).

\bibitem[LW07]{LW07}
Yuk-Kam Lau and Jie Wu, \emph{A large sieve inequality  of %
Elliott-Montgomery-Vaughan type for automorphic forms and two applications}, %
preprint (2007).

\bibitem[Mic02]{mic03}
P.~Michel, \emph{{Analytic number theory and families of automorphic
  $L$-functions}}, %
  in \emph{Automorphic Forms and Applications}, %
  P.~Sarnak \& F.~Shahidi ed., %
  IAS/Park City Mathematics Series, vol.~12, %
  American Mathematical Society, Providence, RI, 2007. %

\bibitem[MV99]{MV99}
H.~L. Montgomery and R.~C. Vaughan, \emph{Extreme values of {D}irichlet
  {$L$}-functions at {$1$}}, Number theory in progress, Vol. 2
  (Zakopane-Ko\'scielisko, 1997), de Gruyter, Berlin, 1999, pp.~1039--1052.
  

\bibitem[Roy01]{roy01}
E.~Royer, \emph{Statistique de la variable al\'eatoire {$L({\rm sym}\sp
  2f,1)$}}, Math. Ann. \textbf{321} (2001), no.~3, 667--687. 

\bibitem[Roy03]{roy03}
E.~Royer, \emph{Interpr\'etation combinatoire des moments n\'egatifs des valeurs
  de fonctions {$L$} au bord de la bande critique}, Ann. Sci. \'Ecole Norm.
  Sup. (4) \textbf{36} (2003), no.~4, 601--620. 

\bibitem[RW03]{raw03}
Dinakar Ramakrishnan and Song Wang, \emph{On the exceptional zeros of
  {R}ankin-{S}elberg {$L$}-functions}, Compositio Math. \textbf{135} (2003),
  no.~2, 211--244. 

\bibitem[RW05]{RW03}
Emmanuel Royer and Jie Wu, \emph{Taille des valeurs de fonctions {$L$} de
  carr\'es sym\'etriques au bord de la bande critique}, Rev. Mat.
  Iberoamericana \textbf{21} (2005), no.~1, 263--312. 

\bibitem[Sar87]{sar87}
Peter Sarnak, \emph{Statistical properties of eigenvalues of the {H}ecke
  operators}, Analytic number theory and Diophantine problems (Stillwater, OK,
  1984), Birkh\"auser Boston, Boston, MA, 1987, pp.~321--331. 

\bibitem[Ser97]{ser97}
Jean-Pierre Serre, \emph{R\'epartition asymptotique des valeurs propres de
  l'op\'erateur de {H}ecke ${T}\sb p$}, J. Amer. Math. Soc. \textbf{10} (1997),
  no.~1, 75--102. 

\bibitem[TW03]{TeW03}
G.~Tenenbaum and J.~Wu, \emph{Moyennes de certaines fonctions multiplicatives
  sur les entiers friables}, J. Reine Angew. Math. \textbf{564} (2003),
  119--166. 

\bibitem[Vil68]{vil68}
N.~Ja. Vilenkin, \emph{Special functions and the theory of group
  representations}, Translated from the Russian by V. N. Singh. Translations of
  Mathematical Monographs, Vol. 22, American Mathematical Society, Providence,
  R. I., 1968. 
\end{thebibliography}
\end{document}